\documentclass[12pt]{amsart}
\usepackage{amssymb}
\input{diagrams.tex}
\diagramstyle[scriptlabels,height=8mm,width=8mm]
\newcommand{\query}[1]%
{\mbox{}\marginpar{\raggedright\hspace{0pt}{\small\em #1}}}%
\theoremstyle{plain}

\newtheorem{thm}{Theorem}[section]
\newtheorem{cor}[thm]{Corollary}
\newtheorem{lem}[thm]{Lemma}
\newtheorem{prop}[thm]{Proposition}

\theoremstyle{definition}
\newtheorem{defi}[thm]{Definition}
\newtheorem{conj}[thm]{Conjecture}
\newtheorem{conv}[thm]{Convention}
\newtheorem{nota}[thm]{Notation}
\newtheorem{rem}[thm]{Remark}
\newtheorem{rems}[thm]{Remarks}
\newtheorem{exa}[thm]{Example}
\newtheorem{exas}[thm]{Examples}
\newtheorem{sit}[thm]{}

\newcommand{\brem}{\begin{rem}}
\newcommand{\brems}{\begin{rems}}
\newcommand{\erem}{\end{rem}}
\newcommand{\erems}{\end{rems}}
\newcommand{\bexa}{\begin{exa}}
\newcommand{\bexas}{\begin{exas}}
\newcommand{\eexa}{\end{exa}}
\newcommand{\eexas}{\end{exas}}
\newcommand{\bdefi}{\begin{defi}}
\newcommand{\edefi}{\end{defi}}
\newcommand{\bcor}{\begin{cor}}
\newcommand{\ecor}{\end{cor}}
\newcommand{\blem}{\begin{lem}}
\newcommand{\elem}{\end{lem}}
\newcommand{\bconv}{\begin{conv}}
\newcommand{\econv}{\end{conv}}
\newcommand{\bconj}{\begin{conj}}
\newcommand{\econj}{\end{conj}}
\newcommand{\bprop}{\begin{prop}}
\newcommand{\eprop}{\end{prop}}
\newcommand{\bthm}{\begin{thm}}
\newcommand{\ethm}{\end{thm}}
\newcommand{\bnota}{\begin{nota}}
\newcommand{\enota}{\end{nota}}
\newcommand{\bsit}{\begin{sit}}
\newcommand{\esit}{\end{sit}}
\newcommand{\be}{\begin{eqnarray}}
\newcommand{\ee}{\end{eqnarray}}

\def\ba{\begin{array}}
\def\ea{\end{array}}
\def\bma{\begin{matrix}}
\def\ema{\end{matrix}}
\def\te{{\tilde e}}
\def\bO{{\bar O}}

\def\lto{\longrightarrow}
\def\hto{\hookrightarrow}

\def\cO{{\mathcal O}}
\def\nor{{\rm norm}}
\newcommand{\Sing}{\operatorname{Sing}}
\newcommand{\Spec}{\operatorname{Spec}}
\newcommand{\Frac}{\operatorname{Frac}}

\newcommand{\Proj}{\operatorname{Proj}}

\newcommand{\Aut}{{\operatorname{Aut}}}
\newcommand{\Cl}{{\operatorname{Cl}}}
\newcommand{\Pic}{{\operatorname{Pic}}}

\renewcommand{\div}{{\operatorname{div}}}
\newcommand{\Stab}{{\operatorname{Stab}}}
\newcommand{\Prin}{{\operatorname{Prin}}}
\newcommand{\Div}{{\operatorname{Div}}}

\newcommand{\A}{{\mathbb A}}

\newcommand{\R}{{\mathbb R}}
\newcommand{\C}{{\mathbb C}}
\newcommand{\Q}{{\mathbb Q}}
\newcommand{\Z}{{\mathbb Z}}
\newcommand{\N}{{\mathbb N}}
\newcommand{\T}{{\mathbb T}}

\title{Normal affine surfaces with $\C^*$-actions}

\author{Hubert Flenner}
\address{Fakult\"at f\"ur Mathematik,
Ruhr Universit\"at Bochum,
Geb.\ NA 2/72,
Universit\"ats\-stra\ss e\ 150,
44780 Bochum, Germany}
\email{Hubert.Flenner@ruhr-uni-bochum.de}

\author{Mikhail Zaidenberg}
\address{Universit\'e
Grenoble I, Institut Fourier, UMR 5582 CNRS-UJF, BP 74,
38402 St.\ Martin
d'H\`eres c\'edex, France}
\email{zaidenbe@ujf-grenoble.fr}

\thanks{ This research started during a visit of the first
author at the Institut Fourier of
the University of Grenoble, and continued during a stay of both of us
at the Max Planck Institute of Mathematics
at Bonn and of the second author at the Ruhr University at Bochum.
The authors thank these institutions for
their support.\\
\mbox{\hspace{11pt}}{\it 2000 Mathematics Subject Classification}:
13A02, 13F15, 14R05, 14L30.\\
\mbox{\hspace{11pt}}{\it Key words}: $\C^*$-action,
graded algebra, affine surface, cyclic quotient singularity}

\begin{document}

\begin{abstract} A classification of affine
surfaces  admitting a $\C^*$-action was given in the work of
Bia\l ynicki-Birula, Fieseler and L. Kaup, Orlik and Wagreich,
Rynes and others. We provide a simple alternative description
of normal quasihomogeneous affine surfaces in terms
of their graded rings as well as by defining equations.
This is based on a generalization of the
Dolgachev-Pinkham-Demazure construction.
\end{abstract}

\maketitle

{\footnotesize \tableofcontents}

\section*{Introduction}

A classification of (normal)
affine surfaces
admitting a $\C^*$-action was given e.g., in \cite{Bia, BiaSo,
OrWa, Pi, BasHa, Ry} and \cite{FiKa1}-\cite{FiKa3}.
Here we obtain a simple alternative description of normal affine
surfaces $V$ with a $\C^*$-action in terms of their
graded coordinate rings as well as by defining equations.
Our approach is based on a generalization of the
Dolgachev-Pinkham-Demazure construction \cite{Do, Pi, De}.
Recall (see \cite{FiKa1}-\cite{FiKa3})
that a $\C^*$-action on a normal affine surface $V$ is called

\smallskip

\noindent {\it elliptic} if it has a unique fixed point
which belongs to the closure of every 1-dimensional orbit,

\smallskip

\noindent {\it parabolic} if the set of its fixed points is
   1-dimensional, and

\smallskip

\noindent {\it hyperbolic}
if $V$ has only a finite number
of fixed points,
and these fixed points are of hyperbolic type,
that is each one of them belongs
to the closure of exactly two 1-dimensional orbits.

\smallskip

In the elliptic case, the complement $V^*$
of the unique fixed point
in $V$ is fibered by the 1-dimensional
orbits over a projective curve $C$.
In the other two cases $V$ is fibered over an affine curve $C$,
and this fibration is invariant under the $\C^*$-action.

Vice versa, given a smooth curve $C$ and
a $\Q$-divisor $D$ on $C$, the
Dolgachev-Pinkham-Demazure construction
provides a
normal affine surface $V=V_{C,D}$ with a $\C^*$-action
such that $C$ is just the algebraic quotient of $V^*$ or of $V$,
respectively.
This surface $V$ is of elliptic type if $C$ is projective
and of parabolic type if $C$ is affine.

We remind this construction in
sections 1 and 2 below.
In section 3 we use it to present any normal affine surface $V$ with a
parabolic $\C^*$-action as a normalization of
the surface $x^d-P(z)y=0$ in $\A_\C^3$ for
a certain $d\in\N$ and
a certain polynomial $P\in\C[t]$
(see Theorem \ref{parepr}).

In section 4 we deal with the hyperbolic case.
We generalize the Dolgachev-Pinkham-Demazure construction
in order to make it work for any hyperbolic $\C^*$-surface.
Instead of one $\Q$-divisor $D$ on a smooth affine curve $C$
as before, it involves now two $\Q$-divisors $D_+$ and $D_-$
on $C$. By our result {\em isomorphism classes of normal
affine hyperbolic $\C^*$-surfaces are in 1-1-correspondence to
equivalence classes of triples
$(C,D_+,D_-)$, where $C$ is a  smooth affine curve and $D_+$, $D_-$
is a pair of $\Q$-divisors on $C$ with $D_++D_-\le 0$; two such
triples $(C,D_+,D_-)$ and $(C,D_+,D_-)$ are considered to be
equivalent if and only if $C\cong C'$ and $D_\pm=D_\pm'\pm D_0$ with
a principal divisor $D_0$;} cf.\ Theorem \ref{noname}. We also
determine the structure of the singularities, the
orbits, the divisor class group and the canonical divisor in
terms of the divisors $D_\pm$, see Theorems \ref{smo}, \ref{prop
type}, \ref{class group} and Corollary \ref{canclass}.

Using our description it is possible to represent
any normal hyperbolic $\C^*$-surface fibered over $C=\A^1_\C$
as the normalization of a surface in $\A_C^4$ given by
$$
x^{dk}-P(t)y=0, \quad x^{ek}z-Q(t)=0\quad\text{and}\quad
y^ez^d-R(t)=0\,,
$$
for certain polynomials $P,Q,R\in\C[t]$ satisfying the
relation $P^eR=Q^d$, where $e$, $d$ are coprime. These
polynomials can be easily computed in terms of the data
$(D_+,\,D_-)$ (see Proposition \ref{algebra A}).
For instance, if the divisor $D_-$ is integral
then this system
reduces to one equation $x^ez-Q(t)=0$ in $\A_\C^3$,
and vice versa. When $k=1$
then it again reduces to one equation
$y^ez^d-R(t)=0$ in $\A_\C^3$.

In Proposition \ref{covering}
we show how the pair $(D_+,\,D_-)$
is transformed when
passing to an equivariant cyclic cover of $V$.
We deduce, in particular, a characterization of
normal hyperbolic $\C^*$-surfaces over $C=\A^1_\C$
with the fractional part of $D_-$ supported at one point,
as normalized cyclic quotients of the surfaces
$x^ez-Q(t)=0$ in $\A_\C^3$.

In the forthcoming paper \cite{FlZa3},
which is actually Part II of the present one,
we will apply these results to give a simple description
of all normal affine $\C^*$-surfaces
equipped in addition by a
$\C^+$-action. In fact, this class consists of all normal affine
surfaces which admit an algebraic group action with an open orbit.

We note that the results of this paper hold {\it m.m.}\ for graded
2-dimensional normal algebras of finite type over a Dedekind
domain.

\section{Generalities on graded rings}

A $\Z$-graded ring $A=\bigoplus_{i\in \Z}A_i$ contains
$A_{\ge 0}=\bigoplus_{i\ge 0} A_i$
and $A_{\le 0}=\bigoplus_{i\le 0} A_i$ as subrings. The following lemma
is ``well known"; in lack of a reference we provide a short argument.

\blem\label{dom0}
If $A=\bigoplus_{i\in Z}A_i$ is a finitely generated $A_0$-algebra, then
so are $A_{\ge 0}$ and $A_{\le 0}$. Moreover, $A$ is normal if and only
if so are both $A_{\ge 0}$ and
$A_{\le 0}$.
\elem

\proof
Reversing the grading interchanges the
subrings $A_{\ge 0}$ and $A_{\le 0}$. Thus it is sufficient to prove
the first part for $A_{\ge 0}$. If $a_{ij}\in A_i$ with $-n\le i
\le n$, $j=1,\ldots,n_i$, is a system of homogeneous generators of
$A$, then $A_{\ge 0}$ is generated (as a module over $A_0$) by
the multiplicatively closed system of monomials
$$a^{k}:=\prod_{i,j} a_{ij}^{k_{ij}}\,,$$ where $
k:=(k_{ij})\in \Z^{N}$ satisfies the inequalities \be\label{raco}
k_{ij}\ge 0,\quad -n\le i \le n,\quad j=1,\ldots,n_i, \qquad
\sum_{i,j} ik_{ij}\ge 0\,.\ee By Gordan's Lemma (see \cite{Od})
the rational polyhedral lattice cone $K\subseteq \Z^{N}$ defined
by (\ref{raco}) is a finitely generated semigroup. Hence the
algebra $A_{\ge 0}$ is generated by a finite system of monomials
$a^k\in A_{\ge 0}$.

Next we show that the subalgebra $A_{\ge 0}$
(and then also $A_{\le 0}$) is normal if so is $A$.
Indeed, the integral closure $(A_{\ge 0})_{\text{norm}}\subseteq
A=A_{\text{norm}}$ is graded.
Take a homogeneous element $x\in (A_{\ge 0})_{\text{norm}}$ of
degree $d:=\deg x$, and let
\be
\label{norma}
x^n+\sum_{i=1}^n b_ix^{n-i}=0,\qquad\text{where}\quad b_i\in
A_{\ge 0}\,,
\ee
be an equation of integral dependence. We may
assume that $b_i$ are also homogeneous,
of degree $\deg b_i=di\ge
0$. Since $\deg b_i \ge 0$ we have $d\ge 0$,
and so $x\in A_{\ge
0}$.

    Conversely, suppose that both $A_{\ge 0}$
and $A_{\le 0}$ are normal. The ring $A\otimes_{A_0}\Frac
(A_0)$ is normal and so is equal to $\Frac(A_0)[u,u^{-1}]$ for a
homogeneous element $u$ of minimal degree $>0$ in $A\otimes_{A_0}\Frac
(A_0)$. Hence $A_\nor$ is contained in this subring of $\Frac A$. If
$f\in A\otimes_{A_0}\Frac (A_0)$ belongs to the normalization $A_\nor$
of $A$ then so does its top homogeneous component.
Thus it is enough to deal with
homogeneous elements.
Let $a$ be such an element satisfying an equation
of integral dependence
(\ref{norma}) over $A$. We may suppose
as above that $b_i\in A_{di}$ $(i=1,\ldots,n)$.
Since $di$ has the same sign as $d:=\deg a$, we have
$a\in (A_{\ge 0})_\nor=A_{\ge 0}$ if $d\ge 0$ and
$a\in (A_{\le 0})_\nor=A_{\le 0}$ if $d\le 0$, respectively.
Anyhow, $a\in A$, whence $A$ is normal, as stated.
\qed\medskip

\bnota\label{gen notations} Let $V=\Spec A$
be a normal affine
surface over $\C$ with an effective $\C^*$-action.
The coordinate ring
$A=\bigoplus_{i\in\Z} A_i$ is then naturally graded
so that $A_i$
is the set of elements of $A$ on which $t\in \C^*$ acts via
$t.f=t^i f$. Thus, $A_0=A^{\C^*}$ is the subalgebra of
invariants, and $A_i\quad(i\neq 0)$ consists of the
quasi-invariants of weight $i$. Up to reversing
the grading we may
assume that $A_+:=\bigoplus_{i>0} A_i\neq 0$.
The subsets $A_+$ and
$A_-:=\bigoplus_{i<0} A_i$ of $A$ are ideals in $A_{\ge 0}$ and
$A_{\le 0}$, respectively.
\enota

The following lemma is well known
(see e.g., \cite{De}, \cite[Lemma 1.5]{FiKa1}).

\blem\label{dom}
{\em (a)} If $A_0\neq \C$ then the
set $M:=\{i\in\Z\,\big|\,A_i\neq 0\}$
coincides either with $\N$ or with $\Z$,  and
$A_i$ is a locally free
$A_0$-module of rank $1$ for all $i\in M$. Moreover, if $u\in \Frac
(A_0)\cdot A_1$
is  a non-zero element then
$$
A\subseteq \Frac (A_0)[u,u^{-1}]\,,
\qquad\text{and even }
\qquad A\subseteq \Frac (A_0)[u]\quad\text{if}\quad M=\N\,.
$$
{\em (b)} In particular, if $A_0\cong \C[t]$
then $A_i$ is a free $A_0$-module of rank 1 for all $i\in M$.
\elem

\proof
(a) The $K_0:=\Frac(A_0)$-algebra $A\otimes_{A_0}K_0$ is a
$1$-dimensional normal graded domain over the field $K_0$. Hence it
is isomorphic to the free polynomial ring $K_0[u]$ or the ring of
Laurent polynomials $K_0[u,\,u^{-1}]$, where $u\in K_0A_d$ and $d>
0$.  As the $\C^*$-action is effective $d=1$, and (a) follows.

(b) follows from \cite[Ch. VII, \S 4, Corollary 2]{Bou}.
\qed\medskip

Lemma \ref{dom}(a) does not hold in general without the assumption
that $A_0\neq \C$ as is seen by
the Pham-Brieskorn surfaces
$V_{p,q,r}:=\{x^p+y^q+z^r=0\}\subseteq\C^3$.

\bsit\label{3cases}
Usually (cf.\ \cite{FiKa1}) one distinguishes between the
following three cases.

\begin{enumerate}
\item[(i)] {\it The elliptic case}: $A_{-}=0,\quad A_0=\C$.
\item[(ii)] {\it The parabolic case}: $A_{-}=0,\quad A_0\neq\C$.
\item[(iii)] {\it The hyperbolic case}: $A_{-}\neq 0$.
\end{enumerate}
\esit

Below we provide more information in each of these cases.

\section{The elliptic case}

In the elliptic case the
$\C^*$-action on $V$ is good. In particular, its fixed point set
$F:=V^{\C^*}$ (which is the zero set of the augmentation ideal
$A_{+}$ of $A$) consists of a unique point called {\it the
vertex} of $V$, and the surface $V$ is smooth outside the vertex.
One considers the smooth projective curve $C:=\mbox{Proj}\,A\cong
V^*/\C^*$, where $V^*:=V\setminus F$, together with the orbit
morphism $\pi: V^*\to C$ (the fibers of $\pi$ are the orbits of
the $\C^*$-action on $V^*$).

A useful class of examples of normal affine surfaces with a good
$\C^*$-action is provided by the affine cones over projective
curves. For an ample divisor $D$ on a smooth
projective curve $C$ the ring
$$
A_{C,D}:=\bigoplus_{k\ge 0} H^0(C,\cO_C(kD) ))\cdot u^k\subseteq
\Frac(C)[u]\,,
$$
where $u$ is an indeterminate, is the
coordinate ring of a normal affine surface $V:=\Spec A_{C,L}$ with a
good $\C^*$-action. This surface $V$ is a cone over $C$ obtained
by blowing down the zero section of the line bundle associated to
$\cO_C(-D)$.

Let furthermore a finite group $G$ act on $V$ freely off the
vertex, and assume that this action commutes with the given good
$\C^*$-action on $V$. Then the quotient $V/G$ is again a normal
affine surface with a good $\C^*$-action. Conversely, the following
result is true.

\bthm \label{cone} (\cite{Do, Pi, De, Ru}) Every normal affine
surface with a good $\C^*$-action appears as the quotient of an
affine cone over a smooth projective curve by a finite group
acting freely off the vertex of the cone. \ethm

Generalizing the construction above, for a smooth
projective curve $C$ and a $\Q$-divisor $D$ on $C$ one considers
the graded ring
$$A_{C,D}:=\bigoplus_{k\ge 0} H^0(C,\cO(\lfloor
kD \rfloor))\cdot u^k\,,$$ where $\lfloor E \rfloor$ denotes the integral
part of a $\Q$-divisor $E$. We have the following result.

\bthm \label{doco} (\cite{Pi}, \cite[Theorem 3.5]{De}) Given a
normal affine surface $V=\Spec A$ with a good $\C^*$-action there
exists a  $\Q$-divisor $D$ on the curve $C=\Proj\,A$ such that
$A\cong A_{C,D}$.\ethm

The affine toric surfaces provide an interesting family
of elliptic $\C^*$-surfaces.

\bexa\label{toric} (\cite{Od, Co}) We remind
that a normal affine toric surface
$V=V_{\sigma}$ is associated to a strictly convex rational
polyhedral cone
$\sigma\subseteq\R^2$. If
$\dim\sigma=0$ or $=1$ then
$V_{\sigma}\cong\C^*\times\C^*$ or
$V_{\sigma}\cong\A^1_\C\times\C^*$, respectively, and so $A^\times \neq
\C^*$. Consequently, these two cannot be elliptic $\C^*$-surfaces.
Otherwise, if
$\dim\sigma=2$ then choosing an appropriate base $e_1,e_2$ of the
lattice one may suppose that
$\sigma$ is the cone $C(e_2, d e_1- ee_2)$, where
$d\ge 1$, $0\le e<d$ and $\gcd(e,d)=1$. We denote
$V_{d,e}:=V_{\sigma}$; then $V_{d,e}=\Spec\, A_{d,e}$, where
$$
A_{d,e}:=\bigoplus_{b\ge 0,\,\, ad-be\ge 0}
\C \cdot x^{a}y^{b}\subseteq \C[x,x^{-1}, y,y^{-1}]
$$
is the semigroup
algebra of the dual cone $\sigma^{\vee}=C\,(e_1,
e e_1+d e_2)$.

The 2-torus
$\T=(\C^*)^2$
acts on $V_{d,e}$ with an open
orbit $V_{d,e}^*:=V_{d,e}\backslash\{\bar 0\}$. Thus one can
introduce on
$V_{d,e}$  a number of elliptic, parabolic as well as
hyperbolic $\C^*$-actions by choosing appropriate 1-parameter
algebraic  subgroups of the torus $\T$.

In \cite{Ri, BenRi1, BenRi2, Co}
one can find a description of minimal sets of generators of the
algebras $A_{d,e}$ as above, as well as defining equations for
the affine varieties $V_{d,e}=\Spec A_{d,e}\hookrightarrow\C^N$.
An explicit presentation of these algebras as in Theorem
\ref{doco} is given in \cite[5.1]{De}.
\eexa

We would like to emphasize the well known relation
between affine toric surfaces and cyclic quotient singularities
(see \cite[5.2]{De} or \cite[Proposition 1.24]{Od}).

\blem\label{LRa}
If $B$ is the normalization of $A:=A_{d,e}$ in the
field $L:=\Frac(A)[u]$ with $u:={\root d \of
x}$, then $B$ is the polynomial ring $B=\C[u, v]$ with $v:=u^ey$.
The Galois group $\langle\zeta\rangle\cong\Z_d$ of $L:\Frac(A)$ acts
on $B$ via the representation, say $G_{d,e}$
$$
\zeta .u=\zeta u,\qquad \zeta .v=\zeta^{e} v\,,
$$
and $A=B^{\Z_d}$. Consequently, there is an isomorphism
$$
V_{d,e}\cong \A^2_\C/G_{d,e}=\A^2_\C/\Z_d\,.
$$
\elem

\proof
For the convenience of the reader we give a short argument.
By definition, $A$ is generated over  $\C$ by the monomials
$$
    x^ay^b \qquad
\text{with} \qquad b\ge 0,\,\,\,ad-be\ge 0\,.
$$
As $x^ay^b=u^{ad-be}v^b$, this shows that $A$ embeds naturally
into $\C[u,v]$ and that even
$A=\C[x,x^{-1},y]\cap \C[u,v].$
In particular $A$ is a normal domain.
Because of $u^d=x\in A$ and
$v^d=x^ey^d\in A$ the ring $B$ is integral over $A$, whence it is
the normalization of $A$.

The second part follows from the first one, since $L$ is a cyclic
extension of $\Frac(A)$ with Galois group $\Z_d$ acting via
$\zeta.u= \zeta u$ and $\zeta. x=x$ for all $x\in A$.
\qed\smallskip

\brem
Assuming that $e>0$ and
letting $\xi:=\zeta^e$ one obtains
$$
(\zeta u, \zeta^{e}v)=(\xi^{e'}u,\xi v)\,,
$$
where $0\le e'<d$ and $ee'\equiv 1\mod d$
(note that for $d=1$ this means $e'=0$). Hence, with
$\tau(u,v):=(v,u)$ the conjugate
$\Z_d$-action
$G_{d,e'}':=\tau^{-1}G_{d,e'}\tau$ on $\A^2_\C$
$$
\xi . (u,v)=(\xi^{e'}u,\xi v)
$$
has the same orbits as $G_{d,e}$
thus providing an isomorphism of toric surfaces
$$
V_{d,e}\cong \A_\C^2/G_{d,e}\cong \A_\C^2/G_{d,e'}'\cong
\A_\C^2/G_{d,e'}\cong V_{d,e'}\,.
$$
Moreover, $V_{d,e}\cong V_{d',e'}$ if and only
if $d=d'$ and either $e=e'$ or $ee'\equiv 1 \mod d$.
\erem

\section{The parabolic case}

In the parabolic case one considers a normal affine surface $V$ with
a $\C^*$-action such that the coordinate ring $A=\bigoplus_{i\ge 0}
A_i$ is positively graded and $A_0$ is a 1-dimensional domain. Thus
$A_0$ corresponds to a smooth affine curve $C=\Spec A_0$, which can
be identified with the algebraic quotient $V//\C^*$ (indeed,
$A_0=A^{\C^*}$ is the ring of invariants of the $\C^*$-action on $A$).
The embedding $A_0\hto A$ corresponds to the quotient morphism
$\pi: V\to C$, and the projection $A\to A_0$ gives an embedding
$\iota: C\hto V$ which provides a retraction of $\pi$ and whose
image is the fixed point set. Every fiber of
$\pi: V\to C$ is the closure of a non-trivial orbit; it contains
a unique fixed point (a {\it source} of this orbit) \cite[Lemma
1.7]{FiKa1}.

A simple example of a parabolic $\C^*$-surface is the cylinder
$C\times\A^1_\C$ over a smooth affine curve $C$, where $\C^*$ acts on the
second factor. More examples can be produced by applying equivariant
affine modifications to $C\times\A^1_\C$ (see \cite[Theorem 1.1]{KaZa1}).
Actually, one
obtains in this way all normal affine surfaces with a parabolic
$\C^*$-action.

\bsit\label{DPD-para}
The Dolgachev-Pinkham-Demazure construction (see Theorem
\ref{doco}) is available also in the parabolic case. Let $C=\Spec
A_0$ be an affine curve over $\C$ with function field $K_0:=\Frac
(A_0)$, and let $D$ be a $\Q$-Cartier divisor on $C$. Similarly as in
the elliptic case we can introduce the algebra
$$
A_0[D]:=A_{C,D}=\bigoplus_{n\ge 0} H^0(C,{\cO}_C(\lfloor nD
\rfloor))\cdot u^n\subseteq K_0[u]\,.
$$
More explicitly, if $f\in K_0$ then \be\label{*} fu^n\in
A:=A_0[D] \Leftrightarrow \div f+ nD\ge 0. \ee We note that the
algebra $A$ is normal (see Corollary \ref{tore}(b) below) 
and finitely generated over $A_0$. Notice also that $u\in A_1$ if
and only if $D\ge 0$. \esit

The following theorem is well known (cf.\ \cite[Theorem
3.5]{De}); for the convenience of the reader
we include a short proof.

\bthm\label{LR}
Let $C=\Spec A_0$ be a normal affine algebraic curve with function
field $K_0:=\Frac(A_0)$. If $A=\bigoplus_{i\ge 0} A_i$ is a
normal finitely generated $A_0$-algebra of dimension 2 then the
following hold.

{\em (a)} $A$ is isomorphic to $A_0[D]$ for some $\Q$-divisor $D$ on
$C$. More precisely,
if  $u\in K_0\cdot A_1$ is a non-zero element and if the divisor
$D$ is defined by the equality
$$
\pi^* D= {\rm div}\, u- \iota(C)\,,
$$
then $A$ and $A_0[D]$  are equal
when considered as subrings of
$K_0[u]$.

{\em (b)}
For two $\Q$-divisors $D$ and $D'$ on $C$,
the rings  $A=A_0[D]$ and $A'=A_0[D']$
are isomorphic as graded $A_0$-algebras
if and only if $D$ and $D'$ are linearly equivalent.
\ethm

\proof
(a) Since $u\in K_0\cdot A_1$ is homogeneous, the divisor
$\mbox{div}\,u$ on the normal surface $V=\Spec A$ is invariant under
the induced $\C^*$-action on $V$, and so we have
$$\mbox{div}\,u=\sum_{i=1}^m p_iF_i + \iota(C)\,$$ with
$p_i\in \Z$, where $F_i=\pi^{-1}(x_i)$ are the fibers of $\pi$
over distinct points $x_i\in C$, $ i=1,\ldots,m$. Letting
$\pi^* x_i= q_iF_i$ with $q_i\in\N\quad (i=1,\ldots,m)$, the
$\Q$-divisor $D:=\sum_{i=1}^m p_i/ q_ix_i$ on $V$
satisfies
$$
\div\,u=\pi^* (D) + \iota(C)\,.
$$
Since $V$ is normal, for a rational function $\varphi\in K_0$ on $C$
the following equivalences  hold:
$$\varphi u^n\in A_n \Leftrightarrow \mbox{div}\,(\varphi u^n)\ge
0 \Leftrightarrow \pi^*\mbox{div}\,\varphi+n\mbox{div}\,u\ge 0
\Leftrightarrow $$
$$\pi^*\mbox{div}\,\varphi + n\pi^*(D)+n\iota(C)\ge 0
     \Leftrightarrow \mbox{div}\,\varphi +nD \ge 0
     \Leftrightarrow \varphi\in H^0(C,{\cO}_C(\lfloor nD \rfloor))\,.$$
Hence $A_n=H^0(C,{\cO}_C(\lfloor nD \rfloor))\cdot u^n$ for all $n\ge
0$, as desired.

(b) Any isomorphism of graded $A_0$-algebras
$$
\varphi: A_0[D]=\bigoplus_{n\ge 0} H^0(C,{\cO}_C(\lfloor nD
\rfloor))\cdot u^n\lto
A_0[D']=\bigoplus_{n\ge 0} H^0(C,{\cO}_C(\lfloor nD'
\rfloor))\cdot u^{\prime n}\,,
$$
extends to an isomorphism of graded $K_0$-algebras
$$
\varphi_{K_0}:K_0[u]\to K_0[u']
$$
and so has the form $u^n\mapsto f^nu^{\prime n}$, $n\ge 0$, for some
non-zero $f\in K_0$. Conversely, such a morphism $\varphi_{K_0}$ maps
$A_0[D]$ isomorphically onto
$A_0[D']$ if and only if
$$
H^0(C,{\cO}_C(\lfloor nD'\rfloor))=
f^n\cdot H^0(C,{\cO}_C(\lfloor nD\rfloor)) \quad \forall n.
$$
As
$$
f^n\cdot H^0(C,{\cO}_C(\lfloor nD\rfloor))=
H^0(C,{\cO}_C(\lfloor nD-n\div f\rfloor)) \,,
$$
the existence of an isomorphism $\varphi$ as above is equivalent to
the existence of an element $f\in K_0$ with $D'=D-\div f$.
\qed\medskip

\bsit\label{frpart}
We denote $\{D\}=D-\lfloor D\rfloor$
the fractional part of a $\Q$-divisor $D$.
Since principal divisors are $\Z$-divisors,
we have $\{D\}=\{D'\}$
as soon as $D \sim D'$.

If $C=\Spec\C[t]=\A^1_\C$ then the converse is also true.
Indeed, any $\Z$-divisor on $\A^1_\C$ is principal, and so
the linear equivalence class of a $\Q$-divisor $D$ on $\A^1_\C$
is uniquely determined by the fractional part $\{D\}$ of $D$.
Thus we obtain the following corollary.
\esit

\bcor\label{frdiv} For every normal
parabolic $\C^*$-surface $V=\Spec A$
with $A=\bigoplus_{n\ge 0} A_n$
and $A_0=\C[t]$, there is
a unique isomorphism  $A\cong A_0[D]$
of graded $A_0$-algebras,
where
$D=0$ or $D=\sum_{i=1}^n \frac{p_i}{q_i} x_i$ with $0<p_i<q_i$,
$\gcd (p_i,q_i)=1\,\,\forall i=1,\ldots,n$
and $x_i\in \A^1_\C$, $x_i\neq x_j$ for $i\neq j$.
\ecor

The next lemma is also well known; in lack of a reference we provide a
short argument.

\blem\label{grge}
Let $D$ be a $\Q$-divisor on a normal
affine variety $S$ and consider the graded ring $A:=\bigoplus_{i\ge
0} A_i$, where $A_i:=H^0(S, \cO_S(\lfloor iD\rfloor))\cdot u^i$.
For $d\in\N$ the following conditions are equivalent.
\begin{enumerate}
\item[(i)]
$dD$ is integral.
\item[(ii)] $A_{d+m}=A_dA_m$ for all $m\ge 0$.
\item[(iii)] The $d$-th Veronese subring
$A^{(d)}:=\bigoplus_{m\ge 0} A_{md}$
is isomorphic to the symmetric algebra $S_{A_0}(A_d)$ i.e.,
$A_{md}=S_{A_0}^mA_d$.
\end{enumerate}
\elem

\proof
Condition $(ii)$ is equivalent to
$$
\cO_S(\lfloor (m+d)D \rfloor)\cong \cO_S(\lfloor mD\rfloor) \otimes
\cO_S(\lfloor dD\rfloor)\qquad \forall m\ge 0\,,
$$
and the latter condition is equivalent to
$$
\lfloor (m+d)D \rfloor=\lfloor mD \rfloor + \lfloor dD \rfloor
\qquad \forall m\ge 0\,.
\leqno{(ii')}
$$
Similarly, $(iii)$ is equivalent to
$$
\lfloor mdD \rfloor=m\lfloor dD \rfloor
\qquad \forall m\ge 0\,.
\leqno{(iii')}
$$
The equivalence of $(i)$,  $(ii')$ and $(iii')$ now follows from
the elementary fact that for a rational number $r=\frac{p}{q}$ and
$d\in \N$  the following conditions are equivalent:
$$
(1)\,\, dr\in\Z\qquad
(2) \,\, \lfloor (m+d)r \rfloor =\lfloor mr \rfloor+\lfloor dr \rfloor
\,\, \forall m\ge 0\qquad
(3)\,\,\lfloor mdr \rfloor= m\lfloor dr \rfloor\,\, \forall m\ge
0.
$$
\qed

\bnota\label{da}
We denote $d(A)$ the smallest positive integer
$d$ satisfying the equivalent conditions of Lemma \ref{grge}.
\enota

\brem\label{alternative} In the situation of Theorem \ref{LR},
one can recover $D$ from the graded ring $A=A_0[D]$ more
algebraically as follows. Consider $d\in \N$ with
$A_dA_i=A_{d+i}$ for all $i\ge 0$ (or, equivalently,
$A_{id}=S^i(A_d)$, see Lemma \ref{grge}) and let $v$ be a
generator of $A_d$ as $A_0$-module; this exists after a suitable
localization of $A_0$. {\em If $u^d=f v$ with $f\in\Frac A_0$,
then $D=\div(f)/d$.} In fact, the ideal $vA$ is equal to $A_{\ge
d}$ and so its zero set has no irreducible components in the
fibers of $\pi$. Thus $\div\,v=d\cdot\iota(C)$ on $V$. Since
$$
\pi^*(D)=\div\, u-\iota(C) \quad \mbox{and}\quad
d\cdot \div\, u= \div\, v+\div\, f\,
$$
as divisors on $V$, we obtain $D=\div(f)/d$.
\erem

A parabolic $\C^*$-surface
$V=\Spec A_0[D]$ has at most cyclic quotient
singularities, as follows from Miyanishi's Theorem (see
\cite[Lemma 1.4.4(1)]{Miy2}).
In the next result (see \cite[Section 5]{De}) we
describe their structure in terms of the divisor $D$.

\bprop\label{tore}
{\em (a)}  If $A_0=\C[t]$ and if $D$ is
supported on the origin in $\Spec A_0=\A^1_\C$ so that
$D=-\frac{e}{d}[0]$ with $\gcd(e,d)=1$, then
$A:=A_0[-\frac{e}{d}[0]]$  is naturally isomorphic
to the semigroup algebra
$$
A_{d,e}=\bigoplus_{b\ge 0,\,\,\,a d-b e\ge 0}
\C\cdot t^{a}v^{b}
$$
graded via $\deg t=0,\,\,\deg v=1$ (cf. Example \ref{toric}).
Consequently, $V:=\Spec A$  is isomorphic to
the toric surface $V_{d,e'}=\Spec A_{d,e'}\cong \A_\C^2/G_{d,e'}$,
where
$e'\equiv e \mod d$ and $0\le e'< d$.

{\em (b)}  If $C=\Spec A_0$ is any normal affine curve over $\C$
and $D$ is a $\Q$-divisor on $C$, then the surface $V=\Spec A_0[D]$
is normal with at most cyclic quotient singularities. More
precisely, if $D(a)=-e/d$ with $\gcd(e,d)=1$ then
$V$ has a quotient singularity of type $(d,e')$ at $\iota(a)$, where
$e'$ is as in (a).
\eprop

\proof
The first part of (a) follows immediately from (\ref{*}) in
\ref{DPD-para}, whereas the second one is a consequence of Lemma
\ref{LRa}.

Tensoring the isomorphism in (a) with
$-\otimes_{\C[t]}\C[[t]]$ we obtain that (b)  holds if $A_0\cong
\C[[t]]$. The general case follows from this by taking completions at
the maximal ideals of $A_0$.
\qed
\medskip

The algebra $A_0[D]$ is
finitely generated over $A_0$, so there exist $f_1,\ldots,f_n\in K_0$
and $m_1,\ldots,m_n\in\N$ such that
$$
A= A_0[f_1u^{m_1},\ldots,f_n u^{m_n}]\subseteq K_0[u]\,.
$$
In the next result we show how to compute $D$ from such a
representation.

\bprop\label{AM}
Let $C=\Spec A_0$ be a smooth affine curve and $K_0:=\Frac A_0$. If
a 2-dimensional subring $B$ of the polynomial ring $K_0[u]$ is
represented as
$$
B=A_0[f_1
u^{m_1},\ldots,f_nu^{m_n}]\subseteq K_0[u]\,,\quad
m_i> 0 \,\,\forall i
$$
with $f_1,\ldots,f_n\in K_0$, then its normalization
$A=B_\nor$ coincides as an $A_0$-subalgebra  of
$K_0[u]$ with $A_0[D]$, where
$$
D:=-\min_{1\le i\le n} \frac{\div\,f_i}{m_i} \,.
$$
\eprop

\proof
By definition of $D$ we have $\div\, f_i+m_iD\ge 0$ so by (\ref{*})
$f_iu^{m_i}\in A_0[D]$ and $B$ is a subring of
$A_0[D]$. As $A_0[D]$ is normal (see Proposition \ref{tore}(b)),
$A$ is also contained in $ A_0[D]$. Let us show that these
subrings coincide.

According to Theorem \ref{LR},
we can represent $A$ as $A=A_0[D']$
with $\pi^*(D')= \div\, u - \iota (C)$. In particular $f_i u^{m_i}\in
A=A_0[D']$, so again by (\ref{*}) $\div f_i +m_iD'\ge 0$ or,
equivalently,
$D'\ge -{1\over m_i} \div\,  f_i.$
Thus $D'\ge D$ and $A_0[D]\subseteq A_0[D']=A$. As we
have already shown the converse inclusion we obtain that
$A=A_0[D]$, as  desired.
\qed\medskip

The following examples
of parabolic $\C^*$-surfaces ruled over $\A^1_\C$
are basic (see Theorem \ref{parepr} below).

\bexa\label{paradp} For a unitary polynomial $P\in\C[t]$
and for an integer $d\ge 1$
we let
$$
B^+_{d,P}:=\C[t,u,v]/(u^d-P(t)v)\cong
\C\left[t,u,{u^d\over P(t)}\right]
$$
graded via
$$
\deg t=0,\qquad \deg u=1,\qquad\deg v=d\,.
$$
The normalization
$$
A^+_{d,P}:=(B^+_{d,P})_\nor
$$
is a positively graded finitely generated $\C$-algebra of
dimension 2 with $A_0=\C[t]$.
By Proposition  \ref{AM} and Corollary \ref{frdiv} we have
$$
A^+_{d,P}\cong A_0[D]\cong A_0[\{D\}],\qquad\text{where}\qquad
D=D(d,P):={\div (P)\over d}\,.
$$
For
$P(t)=\prod_{i=1}^n(t-x_i)^{r_i}$ (where
$x_i\neq x_j$ if $i\neq j$) we obtain
$$
D=\sum_{i=1}^n {r_i\over d} x_i\,,\qquad\text{and}\qquad
\{D\}=\sum_{i=1}^n \left\{{r_i\over d}\right\} x_i\,,
$$
whereas $D=0$ if $P=1$. Replacing $D$ by $\{D\}$
we may suppose that

\medskip

\noindent
$(*)$\quad  $\gcd (d, r_1,\ldots,r_n)=1,\,\,
0<r_i<d\,\,\forall i=1,\ldots,n$, if $d\ge 2$, and
$P=1$ if $d=1 .$

\medskip

\noindent
If two pairs $(d,P)$ and $(\tilde d,\tilde P)$ satisfy $(*)$ and
if $A^+_{d,P}\cong A^+_{\tilde d,\tilde P}$ as graded $A_0$-algebras
then by Corollary \ref{frdiv} we have ${\div (P)\over d}=
{\div (\tilde P)\over \tilde d}$, and so $d=\tilde d$ and $P=\tilde P$.
\eexa

Thus we obtain the following classification result.

\bthm\label{parepr} For every normal affine surface $V=\Spec A$,
where $A=\bigoplus_{i\ge 0} A_i$ with $A_0= \C[t]$,
there is a unique pair $(d,P)$ satisfying condition $(*)$
and an equivariant isomorphism of $A_0$-schemes
$$
\varphi:V\lto V^+_{d,P}:=\Spec A^+_{d,P}\,.
$$
\ethm

\brem\label{3.12}
1.  In the situation of Theorem \ref{parepr} above, the Veronese
subring $A^{(d)}$ is equal to $A_0[v]=\C[t,v]$.
The cyclic group $\Z_d$ acts on $A$ via the $\C^*$-action and
$A^{(d)}$ coincides with the
ring of invariants $A^{\Z_d}$, whereas $A$
is the normalization of $A^{(d)}$ in the fraction field $\Frac (A)$.
Thus the morphism $V\to\A_\C^2=\Spec \C[t,v]$ induced by the
inclusion $\C[t,v]\subseteq A$ represents
$V$ as a cyclic covering of the plane
branched along the curve $u=0$, and $V$ is the normalization of a
surface $\{u^d-P(t)v=0\}$ in $\C^3$.

2.
More generally, let $C=\Spec A_0$ be any smooth affine curve and
let
$A=\bigoplus_{i\ge 0} A_i$ be a normal
2-dimensional $A_0$-algebra of finite type. If $A_1=u\cdot A_0$ and
$A_d=v\cdot A_0$, $d:=d(A)$, for suitable elements $u\in A_1$
and
$v\in A_d$ then
$A$ is the normalization of an algebra
$A_0[u,v]/(u^d-P_+v)$ graded via $\deg u=1,\,\deg v=d$,
for a certain $d\in \N$ and a certain
element $P_+\in A_0$.
\erem

\section{The hyperbolic case}

Let $A=\bigoplus_{i\in\Z}A_i$ be
the coordinate ring of a normal affine surface $V=\Spec A$ with
$\C^*$-action such that $A_+$, $A_-$ are both non-zero. Here again
there is a quotient morphism
$\pi:V\to C=\Spec A_0$ induced by the inclusion $A_0\hookrightarrow
A$. Every fiber of
$\pi$ is either a non-trivial orbit or a union of two
1-dimensional orbits and a hyperbolic fixed point, which is a
source for one of them and a sink for the other one \cite[Lemma
1.7]{FiKa1}. Thus the fixed point set $F$ is finite and contains
$\mbox{Sing}\,V$.

By Lemma \ref{dom0} the proper subalgebras
$A_{\ge 0}$ and $A_{\le 0}$ of $A$ are normal and
finitely generated, and so $V_+:=\Spec A_{\ge 0}$ and $V_-:=\Spec
A_{\le 0}$ are normal affine surfaces with a parabolic
$\C^*$-action birationally dominated by $V$. The natural embeddings
$A_0\hookrightarrow A_{\ge
0}\hookrightarrow A$ and $A_0\hookrightarrow A_{\le
0}\hookrightarrow A$ yield the commutative diagram
\be
\label{diagramm}
\begin{diagram}
V_+ & \lTo^{\sigma_+} & V & \rTo^{\sigma_-}&V_-\\
&\rdTo_{\pi_+} & \dTo^{\pi}&\ldTo_{\pi_-}\\
&& C &&
\end{diagram}
\ee
where $\sigma_\pm$ are equivariant birational morphisms.
Hence $\sigma_\pm$ are equivariant affine
modifications \cite[Theorem
1.1]{KaZa1}. More precisely the following result holds.

\bprop \label{afmo}
$V$ can be obtained from $V_\pm$ by blowing
up a $\C^*$-invariant subscheme and
deleting the proper transform of a
$\C^*$-invariant divisor
$D^\pm$ on $V_\pm$, which contains the fixed point curve
$\iota_\pm(C)\subseteq V_\pm$.
\eprop

\proof

Let us show this for $V_+$,
the proof for $V_{-}$ being similar.
Choose a system of homogeneous generators
$a_1,\ldots, a_n$
of the finitely generated $A_0$-subalgebra
$A_{\le 0}$ and let $f_0\in A_+$ be a non-zero element of degree
$m=-\min_i\deg a_i$. Letting $f_i:=a_i f_0$ for $i=1,\ldots,n$
we obtain
$$
A=A_{\ge 0}\left[{f_1\over f_0},\ldots,{f_n\over f_0}\right]
=
A_{\ge 0}[I/f_0]:=\left\{{x_k\over f_0^k}\mid
x_k\in I^k,\,\,\,k\ge 0\right\}\,,
$$
where $I$ is the graded ideal of $A_{\ge 0}$
generated by $f_0,\ldots,f_n$. Thus
$V=\Spec A$ is obtained by blowing up
$V_+=\Spec A_{\ge 0}$ with center $I$
and deleting the proper transform of
the $\C^*$-invariant divisor
$\div\,f_0$ on $V_+$.
As this divisor contains $\iota_+(C)$, the result follows.
\qed\medskip

For a more precise description of the affine modifications
$\sigma_\pm$ see Remark \ref{preci}.

\bsit\label{DPD} The Dolgachev-Pinkham-Demazure construction is
still available in the hyperbolic case.  In \cite[Theorem 3.5]{De}
it is done under the additional assumption that $A_{-n}\otimes
A_n\to A_0$ is an isomorphism for all $n$. Here we generalize
the construction in order to make it work for any hyperbolic
$\C^*$-surface.

Let $D_+$, $D_-$ be $\Q$-divisors on the smooth
affine curve $C:=\Spec A_0$. For
$n\ge 0$ we consider the $A_0$-submodules
$$
A_{-n}:=H^0(C, \cO_C(\lfloor
nD_-\rfloor))\cdot u^{-n}\quad \mbox{and}\quad
A_n:=H^0(C, \cO_C(\lfloor nD_+\rfloor))\cdot u^n
$$
of $\Frac(A_0)[u, u^{-1}]$, where
$u$ is an indeterminate of degree 1. If $D_++D_-\le 0$ then for $n\ge
m\ge 0$ we have
$$
\lfloor nD_+\rfloor +\lfloor mD_-\rfloor \le \lfloor (n-m)D_+\rfloor,
$$
whence
$A_n\cdot  A_{-m}\subseteq A_{n-m}.$
Similarly, for $0\le n\le m$ we have
$A_n\cdot A_{-m}\subseteq A_{n-m}$. Thus
$$
A:= A_0[D_+,D_-]:=\bigoplus_{n\in \Z} A_n
$$
is a finitely
generated $A_0$-subalgebra of $\Frac(A_0)[u,u^{-1}]$ with $A_{\ge
0}=A_0[D_+]$ and $A_{\le
0}\cong A_0[D_-]$. The grading on $A$ defines a natural hyperbolic
$\C^*$-action on the surface
$V:=\Spec A$. The latter surface is normal
as so are the algebras
$A_0[D_+]$ and $A_0[D_{-}]$
(see Lemma \ref{dom0} and Corollary \ref{tore}(b)).
Conversely, we have the following theorem.\esit

\bthm\label{noname}
If $C=\Spec A_0$ is a smooth affine curve and
$A=\bigoplus_{i\in\Z} A_i$  is a normal graded finitely generated
domain of dimension 2
with $A_{\pm}\neq
0$, then the following hold.

{\em (a)} $A$ is isomorphic to $A_0[D_+,D_-]$, where
$D_+$,
$D_-$ are $\Q$-divisors on $C$ satisfying $D_++D_-\le 0$.
More precisely, if $u\in \Frac(A_0)\cdot A_1$ and if the
divisors $D_+$, $D_-$ on $C$ are defined by
\be\label{divs}
\pi_+^*(D_+)=\div(u)-\iota_+(C)\quad \mbox{and}\quad
\pi_-^*(D_-)=\div (u^{-1})-\iota_-(C),
\ee
where $\pi_\pm$ are as in diagram (\ref{diagramm}) above
and $\iota_\pm:C\hto V_\pm$ are the natural embeddings, then
$D_++D_-\le 0$ and
$A\cong A[D_+,D_-]$.

{\em (b)}  $A_0[D_+,D_-]\cong A_0[D'_+,D'_-]$
as graded $A_0$-algebras if and only if, for a rational
function $\varphi\in\Frac (A_0)$,
one has
$$
D_+'=D_+ +\div\, \varphi\qquad\text{and}\qquad
D_-'=D_- -\div\, \varphi\,.
$$
\ethm

\proof
(a) By Theorem \ref{LR} and its proof we have equalities
$$
A_{\ge 0}= A_0[D_+]\quad \mbox{and}\quad A_{\le 0}= A_0[D_-]
$$
as subalgebras of $\Frac(A_0)[u,u^{-1}]$, whence $A=A_0[D_+,D_-]$.
It remains to show that $D_++D_-\le 0$. Applying in (\ref{divs}) the
functors $\sigma_+^*$ and $\sigma_-^*$ respectively, we obtain
$$
\pi^*(D_+)=\div (u)-\sigma_+^*\iota_+^*(C)\quad \mbox{and}\quad
\pi^*(D_-)=\div (u^{-1})-\sigma_-^*\iota_-^*(C).
$$
Taking the sum of these equalities yields $\pi^*(D_++D_-)=
-(\sigma_+^*\iota_+^*(C)+\sigma_-^*\iota_-^*(C))$, whence
$D_++D_-\le 0$, as required. Finally (b) follows from  Theorem
\ref{LR}(b) and its proof.
\qed\medskip

Consequently, if $A_0=\C[t]$ then $A$
admits a unique presentation $A=A_0[D_+,D_-]$
with $D_+=\{D_+\}$ and $D_++D_-\le 0$.

It follows from Theorem \ref{noname} that outside $|D_+|\cup
|D_-|$, the map $\pi:V\to C$ is  a locally trivial principal
$\C^*$-bundle. More generally, the Dolgachev-Pinkham-Demazure
construction shows the following result (cf.\ \cite{BasHa},
\cite[Proposition 1.11]{FiKa1}).

\bcor\label{cyli}
In all three
cases, outside of a finite subset of the curve $C$ the projection
$\pi:V^*\to C$ and $\pi:V\to C$, respectively,
defines a locally trivial fiber
bundle. This is a principal
$\C^*$-bundle in the elliptic and hyperbolic cases, and a line
bundle in the parabolic case.
\ecor

Note that if $u\in A_1\cup A_{-1}$ is a non-zero
element then its restriction to
a general fiber of $\pi$ gives a fiber coordinate and so a
trivialization over a Zariski open subset of $C$.

\smallskip

\brem\label{rem smooth}
{\em The algebra
$A=A_0[D_+,D_-]$ contains an invertible element of degree $d>0$
if and only if $D_-=-D_+$ and $dD_+$ is a principal divisor on
$C=\Spec A_0$.}
In fact, if $v\in A$ is an invertible
element
of degree $d>0$ then we can write
$$
v=fu^d\in A_d\qquad\text{and}\qquad v^{-1}=f^{-1}u^{-d}\in A_{-d}\,,
$$
where $f\in\Frac (A_0)$ satisfies
$$
\div (f) + dD_+\ge 0\qquad\text{and}\qquad -\div (f) + dD_-\ge 0\,.
$$
Thus $0\ge D_++D_-\ge 0$, whence $D_-=-D_+$. Since $A_d=vA_0$ it
also follows that $dD_+$ is principal. Conversely, if $D_+=-D_-$
and if $dD_+$ is principal, then $vA_0=A_d$ is free
over $A_0$ and $v=fu^d$ with $\div f+dD_+=0$ by Remark
\ref{alternative}. Hence also $\div f^{-1}+dD_-=0$, so
$f^{-1}u^{-d}\in A$ and $v=fu^d$ is a unit in $A$. \erem

The following analogue of Proposition \ref{tore}
holds with a similar proof.

\blem\label{lru} Let $C=\Spec A_0$ be a smooth affine curve with
function field $K_0=\Frac (A_0)$. If a  graded 2-dimensional domain
$B\subseteq K_0[u,u^{-1}]$
is represented as
$$
B=A_0[h_1u^{-n_1},\ldots, h_k u^{-n_k},\,f_1
u^{m_1},\ldots,f_n u^{m_n}]\qquad (\text{where}\quad
n_i,\,m_j>0\,\,\,\forall i,j)
$$
with $h_1,\ldots,h_k,\,f_1,\ldots,f_n\in K_0$
and $B_0=A_0$, then its normalization
$A=B_\nor$ coincides (as a graded $A_0$-subalgebra  of
$K_0[u]$)
with
$A_0[D_+,D_-]$, where
$$
D_-=-\min_{1\le i\le k} {\div \,h_i\over n_i} \qquad\text{and}\qquad
D_+=-\min_{1\le j\le n} {\div\,f_j\over m_j}\,.
$$
\elem

We notice that the assumption $A_0=B_0$
amounts to the inequalities
$$
{\div\,h_i\over n_i} + {\div\,f_j\over m_j}
\ge 0 \quad\forall i,j\,,
$$
which in turn are equivalent to
$D_++D_-\le 0$.

The following lemma provides additional information in the case
that $\lfloor D_\pm\rfloor$ and $d_\pm(A) D_\pm$ are principal
divisors\footnote{or, equivalently, that $A_{\pm 1}$ and $A_{\pm
d_\pm}$ are free $A_0$-modules of rank 1.}.

\blem\label{pmq} Let $A=\bigoplus_{i\in\Z}
A_i=A_0[D_+,D_-]\subseteq \Frac(A_0)[u,u^{-1}]$, and let
$d_\pm=d_\pm (A)$ be the minimal positive integer such that the
divisor $d_\pm D_\pm$ is integral. If $A_{\pm 1}= u_\pm \cdot
A_0$, $A_{\pm d_\pm}= v_\pm \cdot A_0$ and
$$
u_+u_-=Q, \qquad u_\pm^{d_\pm}=P_\pm  v_\pm\,
$$
for some elements $Q,\,P_\pm\in A_0$, then
\be\label{dplmi} D_+={\div\, P_+\over d_+}+D_0\qquad\text{and}\qquad
D_-={\div\, P_-\over d_-}-D_0-\div\,Q\,,
\ee
where $D_0$ is the integral divisor $D_0=\div(u/u_+)$ on
$C=\Spec A_0$. Consequently,
\be\label{qless}
{\div\, P_+\over d_+}+{\div\, P_-\over d_-}\le \div\,Q\,.
\ee
Furthermore, $P_+$ and $P_-$
are uniquely determined by $D_+$ and $D_-$ through
\be\label{frpt}
\{D_+\}={\div\, P_+\over d_+}
\qquad\text{and}\qquad
\{D_-\}={\div\,P_-\over d_-}\,.
\ee
\elem

\proof
We have $u^{d_+}=P_+ \cdot(u/u_+)^{d_+}v_+$ and
$u^{-d_-}=P_-\cdot(u/u_+)^{-d_-}Q^{-d_-}v_-$ and so by Remark
\ref{alternative}
$$
D_+=\frac{\div (P_+ \cdot(u/u_+)^{d_+})}{d_+}=\frac{\div\,P_+}{
d_+}+D_0\,,\quad\mbox{and}
$$
$$
D_-=\frac{\div(P_-\cdot(u/u_+)^{-d_-}Q^{-d_-})}{d_-}
=\frac{\div\,P_-}{ d_-}-D_0-\div \,Q\,.
$$
Now (\ref{qless}) follows
from the inequality $D_++D_-\le 0$. To show (\ref{frpt}), after localizing
$A_0$ we can assume that $P_\pm=S_\pm^{d_\pm} T_\pm$, where
$S_\pm, T_\pm\in A_0$ are elements with
$$
    \div\,S_\pm=\left\lfloor {\div\,P_\pm\over d_\pm} \right\rfloor
\qquad\text{and}\qquad
\div\,T_\pm=\left\{{\div\,P_\pm\over d_\pm}\right\}\,,
$$
respectively.
The relation
$\left({u_\pm/ S_\pm}\right)^{d_\pm}=T_\pm v_\pm$
then shows that $u_\pm/S_\pm$ is integral over $A$ and so by the
normality of $A$ is contained in $A_{\pm 1}$. As $u_\pm$ is a
generator of $A_{\pm 1}$ this forces that $S_\pm\in A_0^\times$ are
units, proving (\ref{frpt}).
\qed\medskip

In many cases the surfaces $V=\Spec A_0[D_+,D_-]$ can be
represented by explicit equations as follows.

\bprop\label{algebra A}
   With the assumptions as in Lemma \ref{pmq} the following hold.

\smallskip

{\em (a)} $A=A_0[D_+,D_-]$ is the normalization of the
$A_0$-algebra \be\label{3eq} B:=A_0[u_-,v_+,v_-]/
\left(u_{-}^{d_-}-v_-P_-,\,v_+^{d'_-}v_-^{d'_+}- P\,,v_+u_-^{d_+}-
Q_+\right)\, \ee graded via $\deg u_-=- 1, \deg v_\pm=\pm d_\pm$,
where $k:=\gcd(d_+,d_-)$, $d'_\pm :=d_\pm/k$ and \be\label{eq4}
P:=\frac{Q^{kd'_+d'_-}}{P_+^{d'_-}P_-^{d'_+}}\in A_0,\qquad
Q_+:=\frac{Q^{d_+}}{P_+}\in A_0\,. \ee

{\em (b)}
$V=\Spec A$ is a cyclic branched covering of degree $k$ of the
normalization of the hypersurface $\{v_+^{d'_-}v_-^{d'_+}-P=0\}$
in
$C\times\A^2_C$.

\smallskip

{\em (c)} If $k=1$ i.e., if $d_+$ and $d_-$ are
coprime and if $v_+$ is not invertible, then $V=\Spec A$ can be
represented as the normalization of a hypersurface $X$ in
$A^3_\C=\Spec \C[s, v_+,v_-]$ with equation
$$
q(s, v_+^{d_-}\cdot v_-^{d_+})=0\,, $$ where $q\in\C[s,t]$ is a
suitable irreducible polynomial. \eprop

\proof (a) First we note that $A$ is
integral over the subring $A_0[v_\pm]$. Indeed, if $w\in A_k$ with
$k\neq 0$ then $w^{d_+}=av_+^k$ if $k>0$ and
$w^{d_-}=av_-^k$ if $k<0$, where $a\in A_0$ (see Lemma
\ref{grge}). Since $A$ and its subring $A_0[u_-,v_\pm]$ have the
same field of fractions, it follows that $A$ is the normalization
of $A_0[u_-,v_\pm]$.

To find the relations between the generators of
$A_0[u_-,v_\pm]$, note that $v_\pm =u_\pm^{d_\pm}/P_\pm$ and so
$$
v_+^{d'_-}v_-^{d'_+}=
\frac{u_+^{d_+d'_-}u_-^{d'_+d_-}}{P_+^{d'_-}P_-^{d'_+}} =
\frac{Q^{kd'_+d'_-}}{P_+^{d'_-}P_-^{d'_+}}=P\in A_0
$$
Similarly
$$
v_+u_-^{d_+}= \frac{u_+^{d_+}u_-^{d_+}}{P_+} =
\frac{Q^{d_+}}{P_+}=Q_+\in A_0\,.
$$
The general fibers of the natural map $\Spec B\to C=\Spec
A_0$ are irreducible, and every fiber is 1-dimensional and in the
closure of the generic fiber. Thus the surface $\Spec B$ is
irreducible, and (a) follows.

(b) $A_0[v_\pm]$ is contained in the Veronese subring
$A^{(k)}$ of $A$ and the fraction fields of both rings coincide.
As $A$ and then also $A^{(k)}$ is integral over $A_0[v_\pm]$
the normalization of $A_0[v_\pm]$ is just $A^{(k)}$. The
cyclic group $ \Z_k$ acts on $A$ via the $\C^*$-action with
invariant ring $A^{(k)}$. Thus $V \to
\Spec A^{(k)}$ is a cyclic branched covering of degree $k$, and
(b) follows.

(c) In case  $k=1$ the algebra $A=A^{(k)}$ is
itself the normalization of the
hypersurface $A_0[v_+,v_-]/(v_+^{d_-}v_-^{d_+}-P)$. Notice that $P$ is
non-constant as $A$ is a domain and, by our assumption, the elements
$v_\pm$ are not invertible. For a general element $s$ of $A_0$ the map
$\varphi=(s,t)$ is a finite morphism of
$C =\Spec A_0$ onto a plane curve $\tilde C\subseteq \A^2_\C$
with an irreducible equation $q(s,t)=0$, where
$t:=P=v_+^{d_-}v_-^{d_+}\in A_0$. This implies (c). 
\qed\medskip

\brems\label{user} 1. It is worthwhile mentioning how to get,
under the assumptions as in (c), a representation $A\cong
A_0[D_+,D_-]$ in terms of $P$ in (\ref{eq4}). Choose $p$, $q\in
\Z$ with $|{d_+\atop d_-}{p\atop q}|=1$ so that $u':=v_+^qv_-^p$
has degree 1. By an easy calculation $u^{\prime d_+}=v_+P^p$ and
$u^{\prime-d_-}=v_-/P^q$, whence by Remark \ref{alternative}
$A\cong A_0[D_+,D_-]$ with
$$
D_+=\frac{p}{d_+}\div\,P,\quad
D_-=-\frac{q}{d_-}\div\,P,\quad\mbox{and}\quad
D_++D_-=-{\div\,P\over d_+d_-}\,.
$$

2. In analogy with (c), any parabolic $\C^*$-surface $V=\Spec A$
with $A=A_0[D]$, where $\lfloor D\rfloor$ and $d(A)D$ are
principal divisors on $C=\Spec A_0$, can be obtained as
the normalization of a surface $u^d-tv=0=q(s,t)$ in
$\A^4_\C=\Spec \C[s,t,u,v]$ graded via $ \deg\, s=\deg\,t=0$, $
\deg\, u=1$, $ \deg\, v=d$, where $q\in\C[s,t]$ is a suitable
irreducible polynomial (see also Remark \ref{3.12}(2)). \erems

The special case $d_+=1$ leads to the following example.

\bexa\label{mei} (Cf. \cite[Example 4.11]{Be2})
For a unitary polynomial $P\in\C[t]$,
we let $A=A_{d, P}=B_{\rm norm}$ be the
normalization of the $\C$-algebra
$$
B=B_{d,P}:=\C[t, u,v]/(u^dv-P(t))
$$
graded via
$\deg t=0,$ $\deg u=1$, $\deg v=-d$
so that the normal affine surface $V:=\Spec A$ is equipped
with a hyperbolic $\C^*$-action.
As $A\cong A_0[u, Pu^{-d}]$ we can
write
$$
A\cong A_0[D_+,D_-],\qquad\text{where}\qquad D_+=0
\quad\text{and}\quad D_-=-\div (P)/d
$$
(see Lemma \ref{lru}). We can recover $P_\pm$ and $Q$ in Lemma
\ref{pmq} as follows. By the construction given there
$P_+=1$ and by (\ref{frpt}) $\{D_-\}=\div(P_-)/d_-$ (note that
$d=d_-$). This gives
\be \label{plusdi}
\div\,P_-=d \left\{ -\frac{\div\,P}{d} \right\}
\qquad\text{and}\qquad
\div\,Q =\frac{\div\,P+\div\,P_-}{d}
\ee
(see (\ref{dplmi})).
In particular,
$$
A_{\ge 0} \cong A_0[u]\cong \C[t,u]\qquad\text{
and}\qquad
A_{\le 0} \cong A^+_{d,P_-}
$$
as graded $A_0$-algebras, where for the second isomorphism we have to
reverse the grading of one of the rings.
\eexa

This discussion provides the following
characterization of the algebras $A_{d,P}$.

\bprop\label{chpr} If $A=A_0[D_+,D_-]$,
where $A_0\cong\C[t]$ and $D_+,\,\,D_-$
are $\Q$-divisors on $\A^1_\C$ with $D_+ +D_-\le 0$,
then the following conditions are equivalent.
\begin{enumerate}
\item[(i)] $D_+$ is integral i.e., $\{D_+\}=0$.
\item[(ii)] $A_{\ge 0} \cong A_0[u]$ as
graded $A_0$-algebras, where $\deg\,u =1$.
\item[(iii)] $A \cong A_{d,P}$ as graded $A_0$-algebras,
where
$D_++D_-=-{\div\,P\over d}$.
\end{enumerate}
\eprop

Next we study the effect of base
change to the Dolgachev-Pinkham-Demazure
representation.

\bprop\label{covering}
Let $C=\Spec A_0$ be an affine curve with function field
$K_0=\Frac(A_0)$ and let
$$
A:=A_0[D_+,D_-]\subseteq K_0[u,u^{-1}]\,,
$$
where $D_\pm$ are $\Q$-divisors on $C$ satisfying
$D_++D_-\le 0$. Let $L$ be the field
$L:=\Frac(A)[\sqrt[d]{tu^b}]$, where $t\in K_0$ and $b,d\in\N$.
If
$A'$ is the normalization of
$A$ in $L$ then the following hold.

1. $A_0'$ is the normalization of $A_0$ in $K_0[s]$ with
$s:=\sqrt[k]{t}$, where
$k:=\gcd(b,d)$.

2. $A'\cong A_0'[D_+',D_-']$ with
$$
D_\pm':=\frac{k }{d}\left(p^*(D_\pm) \pm \beta\div\,s\right) \,,
$$
where $p:C':=\Spec A_0'\to C$ is the projection and $\beta$ is
defined by $\beta b\equiv k \mod d$.
\eprop

\proof We let $b=b'k$ and $d=d'k$.
The normalization $A'$ admits a natural
$\frac{1}{d}$-grading, and the element
$u^*:=\sqrt[d]{tu^b}$  is of degree $b/d=b'/d'$. If we write $k=\beta
b+\delta d$, then the element
$u':= u^{*\beta}u^\delta\in \Frac (A')$
has minimal possible positive degree $1/d'$. Thus
$$
A'\subseteq \Frac(A_0')[u',u^{\prime-1}].
$$
To compute $A_0'$, we note that $u^{*n}u^{-m}$ with $n$, $m\in \N$ has
degree 0 if and only if $nb'/d'=m$. In particular,
$n=n'd'$ is an integer multiple of $d'$. Thus
$K_0':=\Frac A_0'$ is generated over $K_0$ by
$u^{*d'}u^{-b'}=t^{d'/k}$ (i.e., $n'=1$). As $d'$ and
$k$ are coprime, it follows that $s=\sqrt[k]{t}$ also belongs to
$K'_0$ and that this field is actually generated by $s$ over
$K_0$, proving  (1).

After localizing $A_0$ we may assume that there is an element
$v_+\in A$ of degree $d_+=d(A_{\ge 0})$ with $A_{d_+}= v_+A_0$ (see
\ref{da}).  We claim
that then $A'_{sd_+} =v_+^sA'_0$ for all $s\ge 0$. If not, then
for some $s>0$ and some non-unit $x\in A_0'$ the element $v_+^s/x$
belongs to $A'$, so it is integral over $A$ and there is an
equation
$$
\frac{v_+^{sm}}{x^m}+a_{1}\frac{v_+^{s(m-1)}}{x^{m-1}}+\cdots
+a_m=0\,,
$$
where $m\ge 0$ and $a_i\in A_{id_+}$. Thus $a_i=v_+^{si}q_i$ for
some elements $q_i\in A_0$, whence dividing the equation above
by $v_+^{sm}$ we obtain that
$$
\frac{1}{x^m}+q_{1}\frac{1}{x^{m-1}}+\cdots +q_m=0\,.
$$
As $A'_0$ is integrally closed this is only possible if $x\in A'_0$
contradicting the choice of $x$.

Thus $v=v_+$ is an element satisfying the assumptions of Remark
\ref{alternative}, and we compute with it the divisor $D'_+$ as
follows (the calculation for
$D'_-$ is analogous). If we consider the new grading of $A'$ by
assigning to
$u'$ the degree 1, then
$v_+^k$ becomes an element of degree $dd_+$. Moreover, if
$u^{d_+}=P_+v_+$ with $P_+\in K_0$ then by  Remark \ref{alternative}
$D_+={\div (P_+)\over d_+}$. Since
$$
\begin{array}{rcl}
u^{\prime dd_+} &= &(u^{*\beta}u^\delta)^{dd_+}= (tu^b)^{\beta
d_+}u^{\delta d d_+}\\
&=& t^{\beta d_+} u^{d_+(\beta b+\delta d)}= t^{\beta d_+}u^{d_+k}\\
&=& t^{\beta d_+} P_+^kv_+^k\,
\end{array}
$$
we obtain again by Remark \ref{alternative} that on $C'$
$$
D_+'=\frac{\div (t^{\beta d_+} P_+^k)}{dd_+}=\frac{\beta }{d}\div
(t)+
\frac{k}{d}p^*(D_+)\,,
$$
and (2) follows.
\qed\medskip

Let us consider the following important example.

\bexa\label{cops}
With  $A_0:=\C[t]$, suppose that
$D_+=-{e\over d}[0]$ and that $D_-$ is any $\Q$-divisor on
$\A^1_\C=\Spec A_0$.  Applying Proposition \ref{covering}
to $s:=\sqrt[d]{t}$ (i.e.\ $b=0$) we get that the normalization of
$A:=A_0[D_+,D_-]$ in the field $L:=\Frac (A)[s]$ is given by
$$
A'=A_0'[-e[0],D_-']\subseteq \C(s)[u,u^{-1}]\,,
$$
where $A_0'=\C[s]$ and $D_-'=p^*(D_-)$ (as before, $p:\Spec
\C[s]\to\Spec \C[t]$ denotes the projection $s\mapsto s^d$).
The divisor
$D_+'=-e[0]$ being integral we have
$$
A'\cong A_0'[0,D_+'+D_-']\subseteq
\C(s)[\tilde u, \tilde u^{-1}]\,,
$$
where $\tilde u:=s^{e}u$.

More concretely, if $k:=d_-(A)$ and if we choose a unitary polynomial
$Q\in\C[t]$ with
$D_-=-{\div (Q)\over k}$ then $D_+'+D_-'=-{\div (Q(s^d)s^{ke})\over
k}$. By Example \ref{mei} $A'\cong A_{k,P}$ is the
normalization of
\be\label{Bkp}
B_{k,P}=\C[s,\tilde u,v]/\left(\tilde u^kv-P(s)\right)\,,
\quad\mbox{where}\quad P(s):=Q(s^d)s^{ke}\,.
\ee
The field extension $\Frac (A)\subseteq \Frac (A)[s]$
is Galois with  Galois group $\Z_d=\langle\zeta\rangle$, where
$\zeta . s=\zeta s$. Thus
$$
A\cong \left(A_{k,P}\right)^{\Z_d}\,,
$$
and the action of $\zeta$ on $\tilde u=s^eu$ is given
by $\zeta .\tilde u=\zeta^e \tilde u\,.$
Therefore,   the group $\Z_d$ acts on $A_{k,P}$ via
\be\label{actionBkp}
\zeta . s=\zeta s, \quad\zeta . \tilde u=\zeta^e
\tilde u\quad\text{and}\quad
\zeta . v=v\,.
\ee
\eexa

Thus we obtain the following characterization.

\bprop\label{ops} For an algebra $A=A_0[D_+,D_-]$
with $A_0=\C[s]$
the following conditions are equivalent.
\begin{enumerate}
\item[(i)] $\{-D_+\}={e\over d}[0]$,
where $0\le e<d$ and $\gcd (e,d)=1$.
\item[(ii)] $A\cong (A_{k,P})^{\Z_d}$, where $A_{k,P}$ is the
normalization of $B_{k,P}$ in (\ref{Bkp}) and where
$\Z_d=\langle\zeta\rangle$ acts via the formulas in
(\ref{actionBkp}).
\end{enumerate}
\eprop

Like in the parabolic case
$V$ may possess at most cyclic quotient
singularities, as follows from Miyanishi's Theorem (see
\cite[Lemma 1.4.4(1)]{Miy2}). The type of quotient singularities is
determined from the divisors $D_+,D_-$ by the following result. As
before, $C=\Spec A_0$ is a smooth affine curve with function field
$K_0=\Frac A_0$ and $A:=A_0[D_+,D_-]$ with
$\Q$-divisors $D_+$ and $D_-$ on $C$. Denote
$\pi:V=\Spec A\to C$ the canonical projection.

\bthm\label{smo}
\begin{enumerate}
\item[(a)] The set of singular points $\Sing V$
is contained in the fixed point set $F$ which is the zero locus
$F=V(I)$  of the ideal $I:=A_+A+A_-A$ of $A$.
\item[(b)] The map $\pi|F :F \to C$
is injective, and
$\pi(F )=\{a\in C|D_+(a)+D_-(a)<0\}$.
\item[(c)] For a point $a'\in F $ with image
$a:=\pi(a')\in C$ we write
$$
D_+(a)=-{e_+\over m_+}\qquad\text{and}\qquad D_-(a)={e_-\over m_-}
$$
with the convention that
$$
m_+>0,\quad m_-<0,\qquad \gcd(e_+, m_+)=\gcd (e_-, m_-)=1
\qquad\text{and}
$$
$$
m_+=1\quad\text{if}\quad  D_+(a)=0,\qquad
m_-=-1\quad\text{if}\quad D_-(a)=0\,.
$$
Let $p$, $q\in \Z$ with $|{p\atop q}{e_+\atop m_+}|=1$. Then
$a'\in F$ is a quotient singularity of type
$$
(\Delta(a),e),\quad\mbox{where}\quad \Delta(a):=-\left|
\begin{array}{cc}
e_+&e_-\\ m_+&m_-
\end{array}
\right|
\quad \mbox{and} \quad e\equiv
\left|\ba{cc}p&e_-\\q&m_-\ea\right|
\mod \Delta(a)\,.
$$
In particular, $a'\in\Sing V$ if
and only if $\Delta(a)\ne 1$.
\end{enumerate}
\ethm

\proof
As in the proof of Proposition
\ref{tore}(b) we can reduce the statement
to the case that $A_0=\C[t]$ and $|D_+|\cup |D_-|$ is contained in
the origin, so that $D_\pm=\mp e_\pm/m_\pm [0]$.

(a) The set $\Sing V$ is finite and invariant under
the $\C^*$-action. Hence it is contained in the fixed point
set $F$.

(b) The map $A_0\to A/I$ is obviously surjective.
Thus
$$
\pi|F : F = \Spec (A/I)\to C
$$
is a closed embedding. Moreover, $F= \emptyset$
if and only if $1\in I$ if and only if $1=a_+a_-$
for some homogeneous elements of $A$ of opposite degrees, and the
latter happens if and only if $D_++D_-=0$ by Remark \ref{rem
smooth}.

(c) Notice first that the elements
$$
v_+:=t^{e_+}u^{m_+},\,\,\,v_-:=t^{e_-}u^{m_-}\in K_0[u,u^{-1}]
$$
belong to $A$. Indeed, by definition, the ideal $I+tA$ of $A$
(this is just the maximal ideal of the point $a'\in F$) is
generated by the monomials $t^eu^m$ with $(e,m)\in \Z\times \Z$,
where $(e,m)\neq (0,0)$ and
$$
e+mD_+(0)\ge 0 \quad\text{if}\quad m\ge 0,\qquad
e-mD_-(0)\ge 0 \quad\text{if}\quad m\le 0\,.
$$
In other words, $(e,m)$ is an element of the cone  $C:=
\text{C}\,( (e_+,m_+),\, (e_-,m_-))$ generated by the vectors
$(e_\pm,m_\pm)$ in the plane. Hence $A$ is a toric algebra
generated by the semigroup $C\cap \Z^2$, and so is a quotient
$A_{d,e}$ for some $d,\ e\ge 0$ (see Lemma \ref{LRa}). To
determine $d,e$, we must find a basis of $\Z^2$ such that
$(e_+,m_+)$ is one of the basis vectors. This is done as follows.

If we choose $p,q\in\Z$ with $|{p\atop q}{e_+\atop m_+}|=1$, then the
vectors $\te_1:=(e_+,m_+)$ and $\te_2:=(p,q)$ form a basis of
$\Z^2$, and
$$
(e_-,m_-)=\Delta'\te_1+\Delta\te_2, \quad\mbox{where }\ \Delta':=
\left|{p\atop q}{e_-\atop m_-}\right|
\mbox{ and }\ \Delta:=\Delta(0)\,.
$$
As $\te_1$ and $(e_-,m_-)$ form a basis
of the cone $C$, it follows from Lemma \ref{LRa} that $A$ has a quotient
singularity of type $(\Delta,e)$, where $0\le e< d$ and $e\equiv
\left|{p\atop q}{e_-\atop m_-}\right|\mod \Delta$. Note that $\Delta$
and $\Delta'$ are coprime since so are $e_-$ and $m_-$.

The determinant $\Delta$ has always positive sign as
\be\label{delta} D_+(0)+D_-(0)={\Delta\over m_+m_-}\le 0
\qquad\text{and}\qquad m_+>0, \,\,\,m_-<0\,, \ee and so (c)
follows. \qed\medskip

\bcor\label{smth} If $A_{d,P}$
is the normalization of the algebra
$$
B_{d,P}=\C[t, u,v]/(u^dv-P(t))\,,
$$
where $P(t)=\prod_{i=1}^k (t-a_i)^{r_i}$
with $a_i\neq a_j$ for $i\neq j$
(see Example \ref{mei}), then the singular points of
the surface $V_{d,P}=\Spec A_{d,P}$ are
the points $a_i'\in V_{d,P}$ $(1\le i \le k)$,
where $t=a_i,u=v=0$
and $r_i\nmid d$.
\ecor

\proof It was shown in Example \ref{mei}
that $D_-=0$ and $D_+(a_i)=-{r_i\over d}$.
Therefore, $\Delta (a_i)=e_+>1$ if and only if
$r_i\nmid d$, which implies our assertion.
\qed\medskip

In the sequel we use the following notation.

\bdefi\label{deftype}
Let $O=\C^*z$ be the orbit through a point $z\in V\backslash F$.
Following \cite{FiKa1} we say that $O$ is of type $(d,q)$ if $d$
is the order of the stabilizer
$$\Stab_z=\ker (\C^* \to \Aut
\,O)\subseteq\C^*,\qquad\text{so that}\quad
\Stab_z=\langle\zeta\rangle\cong \Z_d\,,$$ and $q\quad (0\le q<
d)$ is determined from the tangent representation of $\Stab_z$
on the tangent plane $T_zV$ via pseudo-reflections
$$\Stab_z\ni\zeta\longmapsto\left(\begin{array} {cc}
1 & 0\\
0 & \zeta^q
\end{array}\right)\,.$$
The orbit $O$ is called {\it principal} if $d=1$ and {\it
exceptional} otherwise (see \cite{FiKa1}-\cite{FiKa3} for a
detailed description of the structure of $V$ near the exceptional
orbits).
\edefi

In the next result we will characterize the orbit types of the
surface $V=\Spec A$ with $A:=A_0[D_+,D_-]$, where $D_+$ and
$D_-$ are $\Q$-divisors on the smooth affine curve $C=\Spec A_0$. Let
$\pi:V\to C$ denote the projection. To examine the orbits over a
point $a\in C$, we write
$$
D_+(a)=- e_+/m_+ \quad\mbox{and}\quad D_-(a)= e_-/m_-
$$
with the conventions as in Theorem \ref{smo}(c). Let $q_+$ be
defined by $0\le q_+ < m_+$ and $q_+ e_+ \equiv -1\mod m_+$, and
similarly $q_-$ by $0\le q_-<-m_-$ and $q_-e_-\equiv 1 \mod
m_-$. With this notation the following result holds.

\bthm\label{prop type}
The exceptional orbits of $V$ are located
over $|D_+|\cup |D_-|$.  The orbits over a given point $a\in
|D_+|\cup |D_-|$ are as follows.

(a) If $D_+(a)+D_-(a)=0$ then $\pi^{*}(a)=m_+O$ consists of one
orbit $O$ of type $(m_+,q_+)$ with multiplicity $m_+$. Moreover,
$O$ appears with coefficient $-e_+$ in $\div\,u$.

(b) If $D_+(a)+D_-(a)<0$ then $\pi^{-1}(a)$ contains two orbits
$O^+$ and $O^-$ of types $(m_+,q_+)$ and $(-m_-,q_-)$,
respectively. Their closures $\bO^\pm$ intersect in the unique
fixed point of the fiber, and $\pi^{*}(a)=m_+\bO^+-m_-\bO^-$.
Moreover, $\bO^\pm$ appears with multiplicity $\mp e_\pm$ in
$\div\, u$. \ethm

\proof With the same reasoning as in the proof of Proposition
\ref{tore}(b) it is sufficient to treat the case where
$A_0=\C[t]$ and $D_\pm$ are supported on $a=0\in \A^1_\C$, i.e.\
$D_\pm=\mp e_\pm/m_\pm[0]$. Note that in this case $m_+=d(A_{\ge
0})$ and $m_-=-d(A_{\le 0})$.

(a) If $D_++D_-=0$, so that $e_+=-e_-=:e$ and $m_+=-m_-=:m$ then
$A$ is the semigroup algebra $\C[C\cap\Z^2]$, where $C$ is the
cone generated over $\R$ by the vectors $\pm(e,m)$ and $(1,0)$. If
we choose $p, q\in \Z$ with $|{p\atop q}{e\atop m}|=1$ then
$$
C\cap \Z^2 =\{(a,b)| (a,b)\in \Z(e,m)+\N(p,q)\}.
$$
Hence $A$ is the algebra of Laurent polynomials \be\label{type1}
A=\C[x,x^{-1},y],\quad\text{where}\quad x:=t^eu^m\in A_m \quad
\mbox{and}\quad  y:=t^p u^q\in A_q\,. \ee Clearly then
\be\label{type2} t=x^{-q}y^m \quad \mbox{and}\quad
u=x^py^{-e}\,.\ee The action of $\C^*$ is given by
$\lambda.x=\lambda^mx$ and $\lambda.y=\lambda^qy$, whence there
is only one orbit $O$ over $t=0$, and it is given by the equation
$y=0$.  By (\ref{type2}) we have
$$
\pi^*(0)=\div\,t=m\cdot O\quad \mbox{and}\quad \div\,u=-e\cdot O.
$$
The stabilizer of any point of $O$ is the group $E_m\subseteq
\C^*$ of $m$-th roots of unity, and the type of the orbit is
$(m,q)=(m_+,q_+)$, as required in (a).

(b) Let now $D_++D_-<0$. Consider a generator $v_\pm=t^{e_\pm}
u^{m_\pm}$ of $A_{m_\pm}$ as $A_0$-module (cf. the proof of
Theorem \ref{smo}(c)). The localization
$A_{v_+}=A[t^{-e_+}u^{-m_+}]$ is the subring $A_0[D_+, -D_-']$ of
$\Frac(A_0)[u,u^{-1}]$ with $D_-':=\max(D_-, -D_+)$ (see\ Lemma
\ref{lru}). As $D_++D_-\le 0$ we have $D_-'=-D_+$, so by (a) the
open subset $\Spec A_{v_+}$ of $V$ contains an orbit $O^+$ of
type $(m_+,q_+)$, and it has multiplicities $m_+$ and $-e_+$ in
$\pi^*(0)$ and $\div\,u$, respectively. Similarly, $\Spec
A_{v_-}$ contains an orbit $O^-$ of type $(-m_-,q_-)$, which has
multiplicities $-m_-$ and $e_-$ in $\pi^*(0)$ and $\div\,u$,
respectively. We have $\div (v_+v_-)=\Delta\cdot \left(\bar
O^++\bar O^-\right)$, where by our assumption
$\Delta=m_+m_-(D_+(0)+D_-(0))>0$ (see (\ref{delta})). Thus the
fiber of $\pi$ over $t=0$ can be given by $v_+\cdot v_-=0$, where
the functions $v_+,\, v_-$ vanish on $\bar O^-$ and $\bar O^+$,
respectively. The intersection $\bar O^+\cap \bar O^-$ is given
by $v_+=v_-=0$, and so is the unique fixed point of the fiber.
\qed

\bexa \label{fxpts} In the example of the algebra $A=A_{d,P}$
treated in Corollary \ref{smth} we have $D_+=0$ and
$D_-=-\div(P)/d=\sum_i -\frac{r_i}{d}[a_i]$ (see Example
\ref{mei}). The exceptional orbits are located over the points
$a_i\in\A^1_\C$, and $\pi^{-1}(a_i)=O^+_i\cup\{a'_i\}\cup O_i^-$,
where $a'_i$ is the unique fixed point of the fiber (located over
the point $(0,0,a_i)$ of $\Spec B_{d,P}\subseteq\C^3$). Applying
Theorem \ref{prop type}, the orbit $O^+_i$ is principal, and if
we write $r_i/d=e_i/m_i$ with $\gcd(e_i,m_i)=1$ then $O^-_i$ is
of type $(m_i, q_i)$, where
$$
q_ie_i\equiv -1\mod m_i\quad\mbox{with}\quad 0\le q_i<m_i\,.
$$
\eexa

\brem\label{preci}
We can now precise the character of the
affine modifications $\sigma_\pm: V \to V_\pm$ as in Proposition
\ref{afmo}. Doing this locally we assume first that $A_0=\C[t]$
and $D_\pm$ is supported on $a=0\in \A^1_\C$. If $D_++D_-=0$ then
$A=A_{\ge 0}[v_+^{-1}]=(A_{\ge 0})_{v_+}$, whence $\sigma_+: V
\to V_+$ is an open embedding and
$V=V_+\backslash V$ is the divisor $\div\,v_+=m_+\iota_+(C)$. In
case $D_++D_-<0$, letting in the proof of Proposition
\ref{afmo}
$f_0:=v_+^{-m_-}$, we obtain that $\sigma_+: V \to V_+$ consists
in blowing up a graded ideal $I\subseteq (t,v_+)$ of the algebra
$A_{\ge 0}$ supported at a fixed point and deleting the proper
transform of the divisor $\div\,v_+=m_+\iota_+(C)$. The
exceptional curve in $V$ is just the orbit $O^-=\{v_+=0\}$.

Globalizing we see that $\sigma_\pm: V \to V_\pm$ blows up a
graded ideal with support at the fixed points
$b_1',\ldots,b_l'\in
\iota_\pm(C)$ over the points $b_i:=\pi_\pm(b_i')\in C$ with
$D_+(b_i)+D_-(b_i)<0$, and deleting the proper transform of the
fixed point curve $\iota_\pm(C)\subseteq V_\pm$. Moreover the
exceptional set of $\sigma_\pm$ is $O_1^\mp\cup\ldots\cup
O_l^\mp$.
\erem

\bsit\label{conve} We let as before $C=\Spec A_0$ be a smooth
affine curve with function field $K_0=\Frac A_0$, and we let
$D_+$, $D_-$ be $\Q$-divisors on $C$. In what follows we compute
the Picard group and the divisor class group of
$A:=A_0[D_+,D_-]$ (see also \cite[Thm.
5.1]{Mori} and \cite[Cor. 1.7]{Wa} for the elliptic case). We
denote by $a_1,\ldots,a_k$ the points in
$C$ for which $D_+(a)=-D_-(a)\ne 0$, and we let $b_1,\ldots,
b_l\in C$ be the points with $D_+(b)+D_-(b)<0$. Let us write
$$
D_\pm(a_i)=\mp\frac{e_i}{m_i},\quad
D_+(b_j)=-\frac{e^+_j}{m_j^+}\quad\mbox{and}\quad
D_-(b_j)=\frac{e^-_j}{m_j^-}\,
$$
with the conventions as in Theorem \ref{smo}.
If $\pi:V:=\Spec
A\to C$ denotes the canonical map then the preimage
$\pi^{-1}(a_i)$ consists of only one orbit $O_i$, and
$\pi^{-1}(b_j)$ consists of two orbit closures
$\bO_j^+\cup\bO_j^-$, so that \be\label{formula fiber}
\pi^*(a_i)=m_iO_i\quad\mbox{and}\quad
\pi^*(b_j)=m_j^+\bO_j^+-m_j^-\bO_j^- \ee
as divisors on $V$, see
Theorem \ref{prop type}.\esit

\bthm\label{class group}
The divisor class group $\Cl\,A$ of $A$ is the group
$$
\pi^*(\Cl\,A_0)\oplus \bigoplus_{i=1}^k\Z[O_i]\oplus
\bigoplus_{j=1}^l\left(\Z[\bO_j^+]\oplus\Z[\bO_j^-] \right)
$$
modulo the relations
$$
\begin{array}{rcl}
\pi^*(a_i)&=&m_i[O_i]\,,\quad i=1,\ldots,k, \\[2pt]
\pi^*(b_j)&=&m_j^+[\bO_j^+]-m_j^-[\bO_j^-]\,,\quad j=1,\ldots l,\\[2pt]
0&=& \sum_{j=1}^k e_i[O_i]+\sum_{j=1}^l
\left(e_j^+[\bO_j^+]-e_j^-[\bO_j^-]\right)\,.
\end{array}
$$
\ethm

\proof
Let $\Div_hA\subseteq \Div\, A$ be the subgroup of all Weil divisors
on $V$ that are homogeneous, i.e.\ finite sums of irreducible divisors
given by homogeneous prime ideals. The homogeneous principal divisors
$\Prin_hA$ form a subgroup of $\Div_hA$, which consists of all
divisors $\div\, f$, where
$f=g/h\in \Frac A$ is a quotient of homogeneous elements.
By \cite[\S 1, Ex.\ 16]{AC}
$$
\Cl\, A\cong\Cl_hA:= \Div_hA/\Prin_hA\,.
$$
The group $\Div_{h}A$ is freely generated by all $\C^*$-invariant
subvarieties of codimension 1 in $V$, that is by all irreducible
components of the fibers of $\pi:V\to C$. If $D_+(a)=D_-(a)=0$
then the fiber over $a$ is the prime divisor $\pi^*(a)$. If
$a=a_i$ for some $i$ then the fiber over $a$ consists of just one
orbit $O_i$ of type $(m_i, q_i)$, and by (\ref{formula fiber})
$\pi^*(a_i)=m_i O_i$ as divisors on $V$. If $a=b_j$ for some $j$
then by (\ref{formula fiber})
$\pi^*(b_j)=m_j^+\bO_j^+-m_j^-\bO_j^-$. Thus the natural map
$\pi^*:\Div\,A_0\to\Div_{h}A$ is injective, and \be\label{formula
div} \Div_hA\cong \frac{ \pi^*(\Div\,A_0)\oplus
\bigoplus_{i=1}^k\Z[O_i]\oplus
\bigoplus_{j=1}^l\left(\Z[\bO_j^+]\oplus\Z[\bO_j^-]\right) }
{(\pi^*(a_i)-m_i[O_i],\,
\pi^*(b_j)-m_j^+[\bO_j^+]+m_j^-[\bO_j^-])}\,\,. \ee The group
$\Prin_hA$ is generated by all divisors $\div(fu^k)=\div\,
f+k\div\,u$, where $f\in K_0^\times$ is non-zero. Dividing out
$\pi^*(\Prin\, A_0)=\pi^*\div(K_0^\times)$ in (\ref{formula div})
gives the group \be\label{formula div1} \frac{
\pi^*(\Cl\,A_0)\oplus \bigoplus_{i=1}^k\Z[O_i]\oplus
\bigoplus_{j=1}^l\left(\Z[\bO_j^+]\oplus\Z[\bO_j^-]\right) }
{(\pi^*(a_i)-m_i[O_i],\,
\pi^*(b_j)-m_i^+[\bO_j^+]+m_j^-[\bO_j^-])}\,\,. \ee By Theorem
\ref{prop type} the divisor of $u$ is given by
$$
\div\, u=
-\sum_{j=1}^k e_i[O_i]+\sum_{j=1}^l
\left(-e_j^+[\bO_j^+]+e_j^-[\bO_j^-]\right)\,.
$$
Hence, taking (\ref{formula div1}) modulo this relation leads to the
divisor class group, as required.
\qed\medskip

\bcor\label{cor class group}
$A$ is factorial if and only if $C\subseteq \A^1_\C$ (i.e.\ $A_0$ is a
localization of $\C[t]$) and one of the following two conditions is
satisfied.
\begin{enumerate}
\item[(i)] $l=0$ and $\gcd(m_i,m_j)=1$ for $1\le i<j\le
k$.
\item[(ii)] $l=1$, $m_i=1$ for all $i$ and $|{e^+\atop
m^+}{e^-\atop m^-}|=\pm 1$, where
$e^\pm:=e^\pm_1$ and $m^\pm:=m^\pm_1$.
\end{enumerate}
\ecor

\proof If $C$ is a curve of genus $g\ge 1$ then the group
$\Cl\,A$ is not finitely generated. Thus assuming that $A$ is
factorial, $C$ is isomorphic to an open subset of $\A^1_\C$. By
Theorem \ref{class group} the group $\Cl\, A$ has then $k+2l$
generators and $k+l+1$ independent relations, whence necessarily
$l\le 1$. In the case $l=1$ the number of generators and the
number of relations are equal, and so the order of $\Cl\,A$ is
the absolute value of the determinant
$$
\left|\begin{matrix}
e^+&e^-&e_1&e_2&\cdots &e_k\\
m^+&m^-&0&0&\cdots&0\\
0&0&m_1&0&\cdots &0\\
0&0&0&m_2&\cdots&0\\
\vdots&\vdots &\vdots&\vdots&\ddots&\vdots\\
0&0&0&0&\cdots& m_k
\end{matrix}\right|
=\left|\begin{matrix}e^+&e^-\\m^+&m^-
\end{matrix}\right|
\cdot m_1\cdot m_2\cdot\ldots\cdot m_k\,.
$$
Thus, if $\Cl\,A=0$ then all the factors of this product are equal
to 1, and we are in case (ii). If $l=0$ then $\Cl\,A$ is the group
$\bigoplus_{i=1}^k\Z_{m_i}\cdot [O_i]$ modulo the relation
$\sum_i e_i[O_i]=0$. As $e_i$ and $m_i$ are coprime, this group is
trivial if and only if (i) holds. Conversely, if (i) or (ii) is
satisfied then the discussion above shows that $\Cl\,A$ is
trivial, finishing the proof. \qed\medskip

Finally, we determine the Picard group and the canonical divisor
of $A$. The local divisor class group at the point $b_j$ is
generated by $\bO^\pm_j$ modulo the relations
$e_j^+\bO_j^+-e_j^-\bO_j^-=0$ and $m_j^+\bO_j^+-m_j^-\bO_j^-=0$.
Since the Picard group $\Pic\, A$ is the kernel of the map of
$\Cl\,A$ into the direct product of all local divisor class
groups, we obtain the following result.

\bcor\label{canclass}
$\Pic\, A$ is the group
$$
\pi^*(\Cl\,A_0)\oplus \bigoplus_{i=1}^k\Z[O_i]\oplus
\bigoplus_{j=1}^l\Z(e^+_j[\bO_j^+]-e_j^-[\bO_j^-])
$$
modulo the relations
$$
\begin{array}{rcl}
\pi^*(a_i)&=&m_i[O_i]\,,\quad i=1,\ldots, k,\\[2pt]
0&=& \sum_{j=1}^k e_i[O_i]+\sum_{j=1}^l
\left(e_j^+[\bO_j^+]-e_j^-[\bO_j^-]\right)\,.
\end{array}
$$
In particular, $\Pic\, A$ vanishes if and only if $C\subseteq
\A^1_\C$ and case (i) in Corollary \ref{cor class group} is
satisfied or $l=1$ and $m_i=1$ for all $1\le i\le k$. \ecor

\bcor
\footnote{Cf. \cite[Thm. 2.8]{Wa} and
\cite[Lemma 2.6]{Mora}.}
The canonical divisor of the surface $V=\Spec A$ is given by
$$
K_V=\pi^*(K_C)+\sum_{j=1}^k (m_i-1)[O_i]+\sum_{j=1}^l
\left((m_j^+-1)[\bO_j^+]+(-m_j^--1)[\bO_j^-]\right)\,.
$$
\ecor

\proof We claim that multiplication by the meromorphic
differential form $du/u$ on $V$ gives an isomorphism
$$
\frac{du}{u}\wedge -: \pi^*(\omega_C)\bigg(\sum_{j=1}^k
(m_i-1)[O_i]+\sum_{j=1}^l
\left((m_j^+-1)[\bO_j^+]+(-m_j^--1)[\bO_j^-]\right)\bigg)
\stackrel{\cong}{\lto} \omega_V\,.
$$
This is a local problem, so with the same arguments as in the
proof of Theorem \ref{prop type} we can reduce to the case that
$A_0\cong \C[t]$ and $D_+=-D_-=-\frac{e}{d}[0]$, where $e$, $m$
are coprime. In this case (\ref{type1}) in the proof of Theorem
\ref{prop type} shows that $A=\C[x,x^{-1},y]$ with $x:=t^eu^m$
and $y:=t^pu^q$, where $p,q$ are integers with $|{p\atop q}{e\atop
m}|=1$. Moreover by (\ref{type2}) $t=x^{-q}y^m$ and
$u=x^py^{-e}$. By an elementary calculation $\frac{du}{u}\wedge
dt= x^{-q-1}y^{m-1}dx\wedge dy$, whence the result follows. \qed


\begin{thebibliography}{mmmm}
\bibitem[BasHa]{BasHa} H. Bass, W. Haboush,
{\em Linearizing certain reductive group
actions}, Trans. Amer. Math. Soc. {\bf 292} (1985), 463-482.

\bibitem[BeRi$_1$]{BenRi1}
K. Behnke, O. Riemenschneider, {\em Diedersingularit\"aten},
Special issue dedicated to the seventieth birthday of Erich
K\"ahler. Abh. Math. Sem. Univ. Hamburg  {\bf 47} (1978), 210-227.

\bibitem[BeRi$_2$]{BenRi2}
K. Behnke, O. Riemenschneider, {\em Quotient surface
singularities and their deformations}, Singularity theory
(Trieste, 1991), 1-54, World Sci. Publishing, River Edge, NJ,
1995.

\bibitem[Be]{Be2} J. Bertin,  {\em  Pinceaux de droites et
automorphismes des surfaces affines}, J. Reine Angew. Math. {\bf
341} (1983), 32--53.

\bibitem[Bi]{Bia} A. Bia\l ynicki-Birula, {\em Some properties of the
decompositions of algebraic varieties determined by actions of a
torus}, Bull. Acad. Polon. Sci. Ser. Sci. Math. Astronom. Phys.
{\bf 24} (1976), 667-674.

\bibitem[BiSo]{BiaSo} A. Bia\l ynicki-Birula, A. J. Sommese, {\em
Quotients by $\C^{*} $ and ${\rm SL}(2,\C)$ actions}, Trans. Amer.
Math. Soc. {\bf 279} (1983), 773-800.

\bibitem[Alg]{Bou} N. Bourbaki, {\em \'El\'ements de
math\'ematique. Alg\`ebre. Chapitres 4 \`a 7}, Lecture Notes in
Mathematics, {\bf 864}, Masson, Paris, 1981.

\bibitem[AC]{AC}  N. Bourbaki, {\em \'El\'ements de
math\'ematique. Alg\`ebre commutative, Chapitre 7: Diviseurs},
Hermann, Paris, 1965.

\bibitem[Co]{Co} D. Cox, {\em Lectures on toric varieties. Lecture 7:
Toric surfaces},
     Grenoble, Institut Fourier, Ecole d'Et\'e 2000 :
     G\'eom\'etrie des vari\'et\'es toriques,
     http://www-fourier.ujf-grenoble.fr/

\bibitem[De]{De} M. Demazure, {\em
Anneaux gradu\'es normaux}, Introduction \`a la th\'eorie des
singularit\'es, II, 35-68, Travaux en Cours {\bf 37}, Hermann,
Paris, 1988.

\bibitem[Do]{Do} I.V. Dolgachev, {\em Automorphic forms and
quasihomogeneous singularities}, Func. Anal. Appl. {\bf 9} (1975),
149-151.

\bibitem[FiKa$_1$]{FiKa1} K.-H. Fieseler, L. Kaup, {\em
On the geometry of affine algebraic $\C\sp *$-surfaces}, Problems
in the theory of surfaces and their classification (Cortona,
1988), 111-140, Sympos. Math., XXXII, Academic Press, London,
1991.

\bibitem[FiKa$_2$]{FiKa2} K.-H. Fieseler, L. Kaup, {\em
Hyperbolic $\C\sp *$-actions on affine algebraic surfaces},
Complex analysis (Wuppertal, 1991), 160-168, Aspects Math., {\bf
E17}, Vieweg, Braunschweig, 1991.

\bibitem[FiKa$_3$]{FiKa3} K.-H. Fieseler, L. Kaup, {\em
Fixed points, exceptional orbits, and homology of affine $\C\sp
*$-surfaces}, Compositio Math. {\bf 78} (1991), 79-115.

\bibitem[FlZa]{FlZa3} H. Flenner, M. Zaidenberg,
{\em Locally nilpotent derivations on affine surfaces with
a $\C^*$-action} (in preparation).

\bibitem[KaZa$_1$]{KaZa1} S. Kaliman, M. Zaidenberg,
{\em Affine modifications and affine varieties with a very
transitive automorphism group}, Transformation Groups, {\bf 4}
(1999), 53-95.

\bibitem[Miy]{Miy2} M. Miyanishi,
{\em Open algebraic surfaces}, CRM Monograph Series {\bf 12},
American Mathematical Society, Providence, RI, 2001.

\bibitem[Mori]{Mori} S. Mori, {\em Graded factorial domains},
Japan.\ J.\ Math.\ {\bf 3} (1977), 223-238.

\bibitem[Mora]{Mora}
M. Morales, {\em Resolution of quasihomogeneous singularities and
plurigenera}, Compositio Math. {\bf 64} (1987), 311-327.

\bibitem[Od]{Od} T. Oda, {\em Convex bodies and algebraic geometry.
An introduction
to the theory of toric varieties}, Ergebnisse der Mathematik und
ihrer Grenzgebiete (3), {\bf 15}, Springer-Verlag, Berlin, 1988.

\bibitem[OrWa]{OrWa} P. Orlik, P. Wagreich, {\em Algebraic
surfaces with $k\sp*$-action}, Acta Math. {\bf 138} (1977), 43-81.

\bibitem[Pi]{Pi} H. Pinkham, {\em Normal surface singularities
with $\C\sp*$ action}, Math. Ann. {\bf 227} (1977), 183-193.

\bibitem[Ri]{Ri} O. Riemenschneider, {\em Die Invarianten der
endlichen Untergruppen von ${\rm GL}(2,\C)$}, Math. Z. {\bf 153}
(1977), 37-50.

\bibitem[Ru]{Ru} B. Runge, {\em Quasihomogeneous singularities}, Math.
Ann. {\bf 281} (1988), 295-313.

\bibitem[Ry]{Ry} J. Rynes,
{\em Nonsingular affine $k\sp *$-surfaces}, Trans. Amer. Math.
Soc. {\bf 332} (1992), 889-921.

\bibitem[Wa]{Wa} Keiichi Watanabe, {\em Some remarks concerning
Demazure's construction of normal graded rings}, Nagoya
    Math.\ J.\ {\bf 83} (1981), 203-211.
\end{thebibliography}
\end{document}